\magnification=1200
%\nopagenumbers

\hsize=11.25cm    
\vsize=18cm       
\parindent=12pt   \parskip=5pt     

\hoffset=.5cm   
\voffset=.8cm   

\pretolerance=500 \tolerance=1000  \brokenpenalty=5000

\catcode`\@=11

\font\eightrm=cmr8         \font\eighti=cmmi8
\font\eightsy=cmsy8        \font\eightbf=cmbx8
\font\eighttt=cmtt8        \font\eightit=cmti8
\font\eightsl=cmsl8        \font\sixrm=cmr6
\font\sixi=cmmi6           \font\sixsy=cmsy6
\font\sixbf=cmbx6

\font\tengoth=eufm10 
\font\eightgoth=eufm8  
\font\sevengoth=eufm7      
\font\sixgoth=eufm6        \font\fivegoth=eufm5

\skewchar\eighti='177 \skewchar\sixi='177
\skewchar\eightsy='60 \skewchar\sixsy='60

\newfam\gothfam           \newfam\bboardfam

\def\tenpoint{
  \textfont0=\tenrm \scriptfont0=\sevenrm \scriptscriptfont0=\fiverm
  \def\rm{\fam\z@\tenrm}
  \textfont1=\teni  \scriptfont1=\seveni  \scriptscriptfont1=\fivei
  \def\oldstyle{\fam\@ne\teni}\let\old=\oldstyle
  \textfont2=\tensy \scriptfont2=\sevensy \scriptscriptfont2=\fivesy
  \textfont\gothfam=\tengoth \scriptfont\gothfam=\sevengoth
  \scriptscriptfont\gothfam=\fivegoth
  \def\goth{\fam\gothfam\tengoth}
  
  \textfont\itfam=\tenit
  \def\it{\fam\itfam\tenit}
  \textfont\slfam=\tensl
  \def\sl{\fam\slfam\tensl}
  \textfont\bffam=\tenbf \scriptfont\bffam=\sevenbf
  \scriptscriptfont\bffam=\fivebf
  \def\bf{\fam\bffam\tenbf}
  \textfont\ttfam=\tentt
  \def\tt{\fam\ttfam\tentt}
  \abovedisplayskip=12pt plus 3pt minus 9pt
  \belowdisplayskip=\abovedisplayskip
  \abovedisplayshortskip=0pt plus 3pt
  \belowdisplayshortskip=4pt plus 3pt 
  \smallskipamount=3pt plus 1pt minus 1pt
  \medskipamount=6pt plus 2pt minus 2pt
  \bigskipamount=12pt plus 4pt minus 4pt
  \normalbaselineskip=12pt
  \setbox\strutbox=\hbox{\vrule height8.5pt depth3.5pt width0pt}
  \let\bigf@nt=\tenrm       \let\smallf@nt=\sevenrm
  \normalbaselines\rm}

\def\eightpoint{
  \textfont0=\eightrm \scriptfont0=\sixrm \scriptscriptfont0=\fiverm
  \def\rm{\fam\z@\eightrm}
  \textfont1=\eighti  \scriptfont1=\sixi  \scriptscriptfont1=\fivei
  \def\oldstyle{\fam\@ne\eighti}\let\old=\oldstyle
  \textfont2=\eightsy \scriptfont2=\sixsy \scriptscriptfont2=\fivesy
  \textfont\gothfam=\eightgoth \scriptfont\gothfam=\sixgoth
  \scriptscriptfont\gothfam=\fivegoth
  \def\goth{\fam\gothfam\eightgoth}
  
  \textfont\itfam=\eightit
  \def\it{\fam\itfam\eightit}
  \textfont\slfam=\eightsl
  \def\sl{\fam\slfam\eightsl}
  \textfont\bffam=\eightbf \scriptfont\bffam=\sixbf
  \scriptscriptfont\bffam=\fivebf
  \def\bf{\fam\bffam\eightbf}
  \textfont\ttfam=\eighttt
  \def\tt{\fam\ttfam\eighttt}
  \abovedisplayskip=9pt plus 3pt minus 9pt
  \belowdisplayskip=\abovedisplayskip
  \abovedisplayshortskip=0pt plus 3pt
  \belowdisplayshortskip=3pt plus 3pt 
  \smallskipamount=2pt plus 1pt minus 1pt
  \medskipamount=4pt plus 2pt minus 1pt
  \bigskipamount=9pt plus 3pt minus 3pt
  \normalbaselineskip=9pt
  \setbox\strutbox=\hbox{\vrule height7pt depth2pt width0pt}
  \let\bigf@nt=\eightrm     \let\smallf@nt=\sixrm
  \normalbaselines\rm}

\tenpoint

\def\pc#1{\bigf@nt#1\smallf@nt}         \def\pd#1 {{\pc#1} }

\catcode`\;=\active
\def;{\relax\ifhmode\ifdim\lastskip>\z@\unskip\fi
\kern\fontdimen2  -1.2 \fontdimen3 \string;}

\catcode`\:=\active
\def:{\relax\ifhmode\ifdim\lastskip>\z@\unskip\fi\penalty\@M\ \fi\string:}

\catcode`\!=\active
\def!{\relax\ifhmode\ifdim\lastskip>\z@
\unskip\fi\kern\fontdimen2  -1.1 \fontdimen3 \string!}

\catcode`\?=\active
\def?{\relax\ifhmode\ifdim\lastskip>\z@
\unskip\fi\kern\fontdimen2  -1.1 \fontdimen3 \string?}

\frenchspacing

\def\raggedbottom{\topskip 10pt plus 36pt\r@ggedbottomtrue}

\def\pointir{\unskip . --- \ignorespaces}

\def\Medbreak{\vskip-\lastskip\medbreak}

\long\def\th#1 #2\enonce#3\endth{
   \Medbreak\noindent
   {\pc#1} {#2\unskip}\pointir{\it #3}\smallskip}

\def\decale#1{\smallbreak\hskip 28pt\llap{#1}\kern 5pt}
\def\decaledecale#1{\smallbreak\hskip 34pt\llap{#1}\kern 5pt}
\def\puce{\smallbreak\hskip 6pt{$\scriptstyle\bullet$}\kern 5pt}

\def\eqalign#1{\null\,\vcenter{\openup\jot\m@th\ialign{
\strut\hfil$\displaystyle{##}$&$\displaystyle{{}##}$\hfil
&&\quad\strut\hfil$\displaystyle{##}$&$\displaystyle{{}##}$\hfil
\crcr#1\crcr}}\,}

\catcode`\@=12

\showboxbreadth=-1  \showboxdepth=-1

\newcount\numerodesection \numerodesection=1
\def\section#1{\bigbreak
 {\bf\number\numerodesection.\ \ #1}\nobreak\medskip
 \advance\numerodesection by1}

\mathcode`A="7041 \mathcode`B="7042 \mathcode`C="7043 \mathcode`D="7044
\mathcode`E="7045 \mathcode`F="7046 \mathcode`G="7047 \mathcode`H="7048
\mathcode`I="7049 \mathcode`J="704A \mathcode`K="704B \mathcode`L="704C
\mathcode`M="704D \mathcode`N="704E \mathcode`O="704F \mathcode`P="7050
\mathcode`Q="7051 \mathcode`R="7052 \mathcode`S="7053 \mathcode`T="7054
\mathcode`U="7055 \mathcode`V="7056 \mathcode`W="7057 \mathcode`X="7058
\mathcode`Y="7059 \mathcode`Z="705A

% handling accented characters in plain TeX :

\def\diagram#1{\def\normalbaselines{\baselineskip=0pt\lineskip=5pt}
\matrix{#1}}

\def\vfl#1#2#3{\llap{$\textstyle #1$}
\left\downarrow\vbox to#3{}\right.\rlap{$\textstyle #2$}}

\def\hfl#1#2#3{\smash{\mathop{\hbox to#3{\rightarrowfill}}\limits
^{\textstyle#1}_{\textstyle#2}}}

\def\ogoth{{\goth o}}

\def\pgoth{{\goth p}}

\def\Q{{\bf Q}}
\def\Qp{\Q_p}

\def\C{{\bf C}}
\def\N{{\bf N}}

\def\Z{{\bf Z}}
\def\Zp{\Z_p}
\def\F{{\bf F}}
\def\Fp{{\F_{\!p}}}

\def\Aut{\mathop{\rm Aut}\nolimits}

\def\Id{\mathop{\rm Id}\nolimits}

\def\Card{\mathop{\rm Card}\nolimits}
\def\Gal{\mathop{\rm Gal}\nolimits}
\def\Ker{\mathop{\rm Ker}\nolimits}

\def\Im{\mathop{\rm Im}\nolimits}

\def\droite#1{\,\hfl{#1}{}{8mm}\,}

\def\to{\rightarrow}

\def\normressym(#1,#2)_#3{\displaystyle\left({#1,#2\over#3}\right)}

\def\mod{\mathop{\rm mod.}\nolimits}
\def\pmod#1{\;(\mod#1)}

\newcount\refno 
\long\def\ref#1:#2<#3>{                                        
\global\advance\refno by1\par\noindent                              
\llap{[{\bf\number\refno}]\ }{#1} \pointir{\it #2} #3\goodbreak }

\def\citer#1(#2){[{\bf\number#1}\if#2\empty\relax\else,\ {#2}\fi]}

\newbox\bibbox
\setbox\bibbox\vbox{\bigbreak
\centerline{{\pc BIBLIOGRAPHIC} {\pc REFERENCES}}

\ref {\pc ARTIN} (E.):
Algebraic numbers and algebraic functions,
<AMS Chelsea Publishing, Providence, 2006. xiv+349 pp.>
\newcount\artin \global\artin=\refno

\ref{\pc DALAWAT} (C.):
Local discriminants, kummerian extensions, and elliptic curves,
<Journal of the Ramanujan Mathematical Society, {\bf 25} (2010) 1,
pp.~25--80. Cf.~arXiv\string:0711.3878v1.>      
\newcount\locdisc \global\locdisc=\refno

\ref{\pc DEL \pc CORSO} (I.) and {\pc DVORNICICH} (R.):
The compositum of wild extensions of local fields of prime degree,
<Monatsh.\ Math.\ {\bf 150} (2007), no.~4, pp.~271--288.>
\newcount\delcorso \global\delcorso=\refno

\ref{\pc FESENKO} (I.) and {\pc VOSTOKOV} (S.):
Local fields and their extensions,
<American Mathematical Society, 2002. xii+345 pp.>
\newcount\fesvost \global\fesvost=\refno

\ref{\pc FONTAINE} (J.-M.):
Groupes de ramification et repr{\'e}sentations d'Artin,
<Ann.\ sci.\ {\'E}cole norm.\ sup.\ (4) {\bf 4}, 1971, pp.~337--392.>
\newcount\fontaine \global\fontaine=\refno

\ref{\pc HASSE} (H.):
Theorie der relativ-zyklischen algebraischen Funktionenk{\"o}rper,
insbesondere bei endlichem Konstantenk{\"o}rper,
<J.\ f.\ d.\ reine und angewandte Math., {\bf 172}, 1934, pp.~37--54.>
\newcount\hasse \global\hasse=\refno

\ref{\pc IWASAWA} (K.):
Local class field theory,
<Oxford University Press, 1986, 155~pp.>
\newcount\iwasawa \global\iwasawa=\refno

\ref{\pc NGUY\~EN-QUANG-\relax\raise1pt\hbox{-}$\!\!$D\~O} (T.):
Filtration de\/ $K^*\!/K^{*p}$ et ramification sauvage,  
<Acta Arith., {\bf 30},  1976, no.~4, pp.~323--340.>
\newcount\nguyen \global\nguyen=\refno

\ref{\pc MARTINET} (J.):
Les discriminants quadratiques et la congruence de Stickelberger,
<S{\'e}m.\ Th{\'e}or.\ Nombres Bordeaux (2) {\bf 1}, 1989,  no.~1,
pp.~197--204.> 
\newcount\martinet \global\martinet=\refno

% \ref{\pc MAUS} (E.):
% Existenz $\pgoth$-adischer Zahlk{\"o}rper zu gegebenem Verzweigungsverhalten,
% <Disseration, Hamburg, 1965.>
% \newcount\maus \global\maus=\refno

\ref{\pc NEUKIRCH} (J.):
Class Field Theory,
<Springer-Verlag, Berlin, 1986, 140~pp.>
\newcount\neukirch \global\neukirch=\refno

\ref{\pc PISOLKAR} (S.):
Absolute norms of\/ $p$-primary numbers,
<J. de Th{\'e}orie de nombres de Bordeaux, {\bf 21} (2009) 3,
pp.~733--740. Cf.~arXiv\string:0807.1174v1.>  
\newcount\supriya \global\supriya=\refno

\ref{\pc SEN} (S.):
Ramification in p-adic Lie extensions,
<Invent.\ Math.\ {\bf 17}, 1972, pp.~44--50.>
\newcount\sen \global\sen=\refno

\ref{\pc SERRE} (J.-P.):
Corps locaux,
<Publications de l'Universit{\'e} de Nancago, No.~{\sevenrm VIII}, Hermann,
Paris, 1968, 245 pp.>
\newcount\serre \global\serre=\refno

\ref{\pc TATE} (J.):
The non-existence of certain Galois extensions of $\Q$ unramified outside
$2$, 
<in Arithmetic geometry, Amer. Math. Soc., Providence, 1994. pp.~153--156.>
\newcount\tate \global\tate=\refno
\ref{\pc WU} (Q.) and {\pc SCHEIDLER} (R.):
The Ramification Groups and Different of a Compositum of Artin-Schreier
Extensions, <Preprint, 2009.>
\newcount\wuscheidler \global\wuscheidler=\refno

} %\bibbox

\centerline{\bf Further remarks on local discriminants} 
\bigskip\bigskip 
\centerline{Chandan Singh Dalawat} 
\centerline{\it Harish-Chandra Research Institute}
\centerline{\it Chhatnag Road, Jhunsi, Allahabad 211019, India} 
\centerline{\it dalawat@gmail.com}

\bigskip\bigskip

{{\bf Abstract}.  Using Kummer theory for a finite extension $K$ of
  $\Qp(\zeta)$ (where $p$ is a prime number and $\zeta$ a primitive $p$-th
  root of~$1$), we compute the ramification filtration and the discriminant of
  an arbitrary elementary abelian $p$-extension of $K$.  We also develop the
  analogous Artin-Schreier theory for finite extensions of $\Fp(\!(\pi)\!)$,
  where $\pi$ is transcendental, and derive similar results for their
  elementary abelian $p$-extensions.  \footnote{}{Keywords~: Local fields,
    Kummer theory, ramification filtration, Artin-Schreier theory, elementary
    abelian $p$-extensions, discriminants.}}

\bigskip\bigskip\bigskip

{\bf 1.  Review of local Kummer theory.}  Before stating the main results of
this paper in \S2, let us briefly review Kummer theory of local fields as
expounded in \citer\locdisc().  This review is justified on the grounds that
\S4 consists largely of applications of this theory, and that \S6 gives
the analogous Artin-Schreier theory for local function fields, which is
summarised in \S5.

We fix a prime number $p$ and denote by $\zeta$ a primitive $p$-th root
of~$1$.  Let $K|\Q_p(\zeta)$ be a finite extension, $k$ the residue field,
$e_1$ the ramification index, and $f=[k:\Fp]$ the residual degree~; the
ramification index of $K|\Qp$ is $e=(p-1)e_1$ (and the residual degree is
$f$).

The filtration $(U_n)_{n>0}$ on $K^\times$ by units of various levels induces
a filtration on the $\Fp$-space $\overline{K^\times}=K^\times\!/K^{\times p}$
denoted by $(\bar U_n)_{n>0}$~; we have $\bar U_{pe_1+1}=\{1\}$, and the
codimension at each step is given by
$$
\{1\}
\subset_1\bar U_{pe_1}
\subset_f\bar U_{pe_1-1}
\subset_f\cdots
\subset_f\bar U_{pi+1}
=\bar U_{pi}
\subset_f\cdots
\subset_f\bar U_1 
\subset_1\overline{K^\times}.
$$
Here, $i$ is any integer in the interval $[1,e_1[$ (which is empty when
$e_1=1$), and an inclusion $E\subset_rE'$ means that $E$ is a codimension-$r$
subspace of $E'$.

We have $\overline{\ogoth^\times}=\bar U_1$, and the valuation gives
$\overline{K^\times}\!/\bar U_1=\Z/p\Z$. Moreover, the choice of $\zeta$ leads
to an isomorphism $\bar U_{pe_1}\to\Fp$ sending $a$ to $S_{k|\Fp}(\hat c)$,
where $\hat c$ is the image of $c=(1-b)/p(1-\zeta)$ in $k/\wp(k)$ and $b\in
U_{pe_1}$ represents $a$~; the isomorphism $k/\wp(k)\to\F_p$ is induced by the
trace map $S_{k|\Fp}$. (Question~: How does the isomorphism $\bar
U_{pe_1}\to\Fp$ change when we replace $\zeta$ by $\zeta^a$ for some
$a\in\F_p^\times$~?) 

The unramified degree-$p$ extension of $K$ is $K(\!\root p\of{U_{pe_1}})$.
For an $\F_p$-line\/ $D\neq\bar U_{pe_1}$ in $\overline{K^\times}$ such that
$D\subset\bar U_m$ but $D\not\subset\bar U_{m+1}$ (with the convention that
$\bar U_0=\overline{K^\times}$), the unique ramification break of the (cyclic,
degree-$p$) extension $K(\!\root p\of D)$ occurs at $pe_1-m$.  This integer is
thus prime to~$p$, unless $m=0$.

Let $M=K(\!\root p\of{K^\times})$ be the maximal elementary abelian
$p$-extension of $K$, and $G=\Gal(M|K)$, endowed with the ramification
filtration $(G^u)_{u\in[-1,+\infty[}$ in the upper numbering.  It follows from
the foregoing (see Part~IX of \citer\locdisc(), cf.~\citer\delcorso()), that
$G^u=G^1$ for $u\in\;]-1,1]$ and that, for $u\in[1,pe_1+1]$, we have the
``\thinspace orthogonality relation\thinspace''
$(G^u)^\perp=\bar U_{pe_1-\lceil u\rceil+1}$, where the orthogonal is taken
with respect to the Kummer pairing $\overline{K^\times}\times G\to{}_p\mu$.
In particular, $G^u=\{\Id_M\}$ for $u>pe_1$.

This shows that the upper ramification breaks of $M|K$ occur precisely at
$-1$, at the $e$ integers in $[1,pe_1[$ which are prime to~$p$, and at $pe_1$.

\bigbreak

{\bf 2.  The main results.}  We begin with a brief account (\S3) of the
congruence satisfied by the absolute norm of a $p$-primary unit, as worked out
by S.~Pisolkar \citer\supriya(). Although not a direct consequence of the
orthogonality relation, her result was inspired by these ideas.

Next, we provide some applications (\S4) of the orthogonality relation.  These
include the computation of the discriminant of elementary abelian
$p$-extensions, the existence of such extensions with given ramification
breaks, their possible degrees and their total number.  It would be possible
to derive some of these results from local class field theory, but our
approach via Kummer theory has the advantage of being more elementary.

Guided by the fecund analogy with local function fields, we then look for an
orthogonality relation for the Artin-Schreier pairing.  The results are proved
in \S6 and summarised in \S5, which should be compared with the review of
Kummer theory in \S1.  

For an opinionated presentation of most of the background, see
\citer\locdisc(), which is freely available online.  Our main theme here is
the compatibility of the Kummer (resp.~the Artin-Schreier) pairing with the
filtration on the multiplicative group $K^\times$ (resp.~ the additive group
$K$) on the one hand and the ramification filtration --- in the upper
numbering --- of $G=\Gal(M|K)$ on the other, where $M$ is the maximal
elementary abelian $p$-extension of a local number field $K$ containing a
primitive $p$-th root of~$1$ (resp.~a local function field $K$) of residual
characteristic~$p$.

This compatibility is expressed by the orthogonality relation as recalled in
\S1 (resp.~proved in \S6).  In essence, it says that $M^{G^n}$ is the same as
$K(\!\root p\of U_m)$ (resp.~$K(\wp^{-1}(\pgoth^{m}))$) whenever $m+n=pe_1+1$
(resp.~$m+n=1$).  The proof is entirely elementary and purely conceptual.

The orthogonality relation has many consequences, some of which we have
discussed above.  It implies, without invoking Hasse-Arf, that the
ramification breaks of $M|K$ occur at integers.  It allows us to compute the
discriminant of any elementary abelian $p$-extension of local fields, without
invoking class field theory and the {\it
  F{\"u}hrerdiskriminantenproduktformel\/}.

When the two approaches are combined, one can compute the norm group of the
extension $K(\!\root p\of U_m)$ (resp.~$K(\wp^{-1}(\pgoth^{m}))$), as
explained in~\S6.

\bigbreak

{\bf 3.  Absolute norms of $p$-primary numbers.}  At the {\it Journ{\'e}es
  arithm{\'e}tiques\/} in Exeter (1980), J.~Martinet generalised the
congruence $D\equiv0,1\pmod4$ for the absolute discriminant of a number field
to a congruence for the absolute norm of the relative discriminant of an
extension of number fields.  One of his results \citer\martinet(p.~198)
suggested the following local version~: if $K|\Q_2$ is a finite extension
containing the $2^m$-th roots of~$1$, and if $L|K$ is a finite unramified
extension, then the discriminant $d_{L|K}\in\ogoth_K^\times$ of any
$\ogoth_K$-basis of $\ogoth_L$ satisfies
$N_{K|\Q_2}(d_{L|K})\equiv1\pmod{2^{m+1}}$.

More generally, it suggested that the absolute norm of any $p$-primary unit in
a finite extension $K$ of $\Q_p$ containing a primitive $p^m$-th root of~$1$
(for some $m>0$) should be $\equiv1\pmod{p^{m+1}}$, where a unit $\alpha$ is
called ``\thinspace$p$-{\it primary}\thinspace'' if the extension $K(\!\root
p\of\alpha)$ is unramified.  This has been verified by S.~Pisolkar~; we
present a variant of her proof.

\th THEOREM {1} (\citer\supriya())
\enonce
Let\/ $K|\Q_p$ be a finite extension for which\/ $K^\times$ has an element of
order\/ $p^m$ ($m>0$), and let\/ $\alpha\in\ogoth_K^\times$ be a unit such
that the extension\/ $K(\!\root p\of\alpha)|K$ is unramified.  Then\/
$N_{K|\Q_p}(\alpha)\equiv1\pmod{p^{m+1}}$. % $N_{K|\Q_p}(u)\in1+p^{m+1}\Z_p$. 
\endth
The proof has four ingredients.  First, we may assume that $K=\Q_p(\xi_m)$,
where $\xi_m\in K^\times$ has order $p^m$.  Assuming this, $\alpha$ can be
written $\alpha=\beta\gamma^p$, where $\beta\in U_{p^m}$ and
$\gamma\in\ogoth_K^\times$.  Thirdly, $N_{K|\Q_p}(\gamma)\equiv1\pmod{p^m}$.
Finally, \hbox{$N_{K|\Q_p}(\beta)\equiv1\pmod{p^{m+1}}$}.  Granting these, the
theorem follows because \hbox{$N_{K|\Q_p}(\gamma)^p\equiv1\pmod{p^{m+1}}$}
(cf.~\citer\locdisc(prop.~{27})).

{\it Reduction to the case\/ $K=\Q_p(\xi_m)$}.  We shall prove that if $F$ is
a finite extension of $\Q_p(\xi_m)$, $E$ is a finite extension of $F$, and
$a\in\ogoth_E^\times$ is a $p$-primary unit of $E$, then
$N_{E|F}(a)\in\ogoth_F^\times$ is a $p$-primary unit of $F$.  We may assume
that $E|F$ is either unramified or totally ramified.

If $E|F$ is unramified, then (\citer\locdisc(prop.~{37})) the image of
$N_{E|F}(a)$ in $F^\times\!/F^{\times p}$ lies in the $\F_p$-line which gives
the unramified degree-$p$ extension of $F$ (\citer\locdisc(prop.~{17})).

Suppose that $E|F$ is totally ramified, and let $E'|E$ be an unramified
extension such that $a=b^p$ for some $b\in E^{\prime\times}$.  There exists an
unramified extension $F'|F$ such that $E'=EF'$~; it suffices to show that
$N_{E|F}(a)\in F^{\prime\times p}$.  Indeed,
$N_{E|F}(a)=N_{E'|F'}(a)=N_{E'|F'}(b)^p$~; the first equality holds because
$N_{E|F}(a)$ (resp.~$N_{E'|F'}(a)$) is the determinant of the
multiplication-by-$a$ automorphism of the $F$-space $E$ (resp.~of the
$F'$-space $E'=E\otimes_F F'$).

\medbreak 

From now on, we may and do assume that $K=\Q_p(\xi_m)$ and denote by $(U_n)_n$
(resp.~$(V_n)_n$) the filtration on $K^\times$ (resp.~$\Q_p^\times$).

{\it One may write\/ $\alpha=\beta\gamma^p$ ($\beta\in U_{p^m}$,
  $\gamma\in\ogoth_K^\times$)}.  This follows from the fact that
$\bar\alpha\in\bar U_{pe_1}$ (\citer\locdisc(prop.~{17})) and the fact that
$e=\varphi(p^m)=p^{m-1}(p-1)$ (\citer\locdisc(prop.~{23})).

\medbreak

{\it Proof that\/ $N_{K|\Q_p}(\ogoth_K^\times)=V_m$.}  We borrow the argument
from \citer\artin(p.~208)~; see also \citer\neukirch(p.~45).  As $K|\Q_p$ is a
totally ramified abelian extension of degree~$\varphi(p^m)$
(cf.~\citer\locdisc(prop.~{23})), the subgroup
$N_{K|\Q_p}(\ogoth_K^\times)\subset\Z_p^\times$ has index~$\varphi(p^m)$
\citer\serre(p.~196).  So has the subgroup $V_m$.  It thus suffices to prove
the inclusion $V_m\subset N_{K|\Q_p}(\ogoth_K^\times)$ to show their equality.

If $p\neq2$, the raising-to-the-exponent-$p$ map $(\ )^p$ is an isomorphism
$V_r\rightarrow V_{r+1}$ for every $r>0$ (\citer\locdisc(prop.~{30})), so
$V_m=V_1^{\varphi(p^m)}$ and the inclusion $V_m\subset
N_{K|\Q_p}(\ogoth_K^\times)$ is clear.

When $p=2$, we may assume that $m>1$ (since $N_{\Q_2|\Q_2}(\Z_2^\times)=V_1$).
Squaring is an isomorphism $V_r\rightarrow V_{r+1}$ for $r>1$
(\citer\locdisc(prop. 30)), so $V_2^{2^{m-2}}=V_m$.  But notice that
$$
V_2=V_3\cup5V_3=V_2^2\cup5V_2^2.
$$
Raising to the exponent $2^{m-2}$ gives
$V_m=V_2^{\varphi({2^m})}\cup5^{2^{m-2}}V_2^{\varphi({2^m})}$.  Clearly,
$V_2^{\varphi({2^m})}\subset N_{K|\Q_2}(\ogoth_K^\times)$.  To get the
inclusion $V_m\subset N_{K|\Q_2}(\ogoth_K^\times)$, it remains to show that
$5^{2^{m-2}}\in N_{K|\Q_2}(\ogoth_K^\times)$.  Indeed, putting
$i=\xi_m^{2^{m-2}}$, so that $i^2=-1$, we have
$$
N_{K|\Q_2}(2+i)=N_{\Q_2(i)|\Q_2}(2+i)^{2^{m-2}}=5^{2^{m-2}}.
$$

(More generally, let $F$ be a finite extension of $\Qp$, let $\pi$ be a
uniformiser of $F$, and let $m$ be a positive integer.  There is a unique
abelian extension $E|F$ such that $\pi\in N_{E|F}(E^\times)$ and
$N_{E|F}(\ogoth_E^\times)=U_{m,F}$ \citer\neukirch(p.~45).  When $F=\Qp$ and
$\pi=p$, then $E=\Qp(\xi_m)$, in view of the fact that $N_{E|\Qp}(1-\xi_m)=p$, 
and, as we have just seen, $N_{E|\Q_p}(\ogoth_E^\times)=V_m$.)

\medbreak

{\it Proof that\/ $N_{K|\Q_p}(U_{p^m})\subset V_{m+1}$.}  We adopt the
notation $\phi=\phi_{K|\Qp}$ and $\psi=\phi^{-1}$ (``Hasse-Herbrand'') for the
piecewise-linear increasing bijections of $[-1,+\infty[$ relative to the
(galoisian ) extension $K|\Q_p$ \citer\serre(p.~73) and use the fact that
$N_{K|\Q_p}(U_{\psi(m)+1})\subset V_{m+1}$ \citer\serre(p.~91).  It thus
suffices to show that $\psi(m)=p^m-1$, or, equivalently, that $\phi(p^m-1)=m$.

The upper ramification subgroups of $G=(\Z/p^m\Z)^\times=\Gal(K|\Q_p)$ are
given by $G^w=\Ker((\Z/p^m\Z)^\times\rightarrow(\Z/p^w\Z)^\times)$ for
$w\in[0,m]$ \citer\serre(p.~79)~; notice that $G^1=G^0$ when $p=2$ because
$\F_2^\times$ is trivial.  The orders are $g^0=p^{m-1}(p-1)$ and $g^w=p^{m-w}$
for $w\in[1,m]$.  The lower indexing is given by
$$
\vbox{\halign{&\hfil$#$\hfil\quad\cr
u\in&\{0\}&[1,p[&[p,p^2[&\cdots&[p^{m-1},+\infty[\cr
\noalign{\vskip-5pt}
\multispan6\hrulefill.\cr
G_u=&G^0&G^1&G^2&\cdots&G^m\cr
}}
$$
(Incidentally, this gives the valuation of the absolute discriminant of $K$ as
$(p^m-p^{m-1}-1)+(p-1)(p^{m-1}-1)+\cdots+(p^{m-1}-p^{m-2})(p-1)$, which equals
$m\varphi(p^m)-p^{m-1}$ and vanishes precisely when $m=1$, $p=2$.  Notice that
when $m>1$, $K$ is a ramified degree-$p$ kummerian extension of
$F=\Q_p(\xi_m^p)$~; its unique ramification break occurs at $p^{m-1}-1$
(cf.~\citer\locdisc(), prop.~{26} and prop.~{60}).  This can also be seen
directly by remarking that $\ogoth_K=\ogoth_F[\xi_m]$ (\citer\locdisc(),
prop.~{23}), that, $\sigma$ being a generator of $\Gal(K|F)$, one has
$\sigma(\xi_m)=\zeta\xi_m$ for some order-$p$ element $\zeta\in K^\times$, and
that the valuation of $1-\zeta$ in $K$ is $p^{m-1}$.  Therefore
$$
v(\sigma(\xi_m)-\xi_m)=v(\zeta\xi_m-\xi_m)=v(1-\zeta)=p^{m-1}
$$
and the ramification break occurs at $p^{m-1}-1$ \citer\serre(p.~61).
Equivalently (cf.~\citer\locdisc()~prop.~{60}) in view of the fact that
$e_1=p^{m-1}$, the image of $\xi_m$ in $K^\times\!/K^{\times p}$ is in $\bar
U_1$ but not in $\bar U_2$.)

Let $g_u$ be the order of $G_u$.  Recall that $g_0\phi(n)=g_1+\cdots+g_n$
\citer\serre(p.~73).  We thus have, for every integer $n\in[1,m]$,
$$
\eqalign{
g_0.\phi(p^n-1)&=
(g_1+\cdots+g_{p-1})
+\cdots
+(g_{p^{n-1}}+\cdots+g_{p^n-1})\cr
&=(p-1).g^1+\cdots+(p^n-p^{n-1}).g^n\cr
&=n.(p^m-p^{m-1})=n.g_0\cr
}$$
and hence $\phi(p^n-1)=n$~; in particular $\phi(p^m-1)=m$, as was to be
proved.  (The same result can also be derived directly from the integral
expression for $\psi$ (\citer\serre(p.~74) or $(3)$ below), which gives
$\psi(1)=p-1$, and, recursively, $\psi(n+1)=\psi(n)+p^{n}(p-1)$ for
$n\in[1,m[$.)  

This completes the proof of Pisolkar's result saying that the absolute norm of
a $p$-primary unit in a finite extension $K|\Q_p$ containing a primitive
$p^m$-th root of~$1$ for some $m>0$ is $\equiv1\pmod{p^{m+1}}$.  The case
$p=2$ of th.~{1} provides a purely local proof of Martinet's generalisation
\citer\martinet() of Stickelberger's congruence.

\bigbreak

{\bf 4.  Elementary abelian $p$-extensions.}  Let $K$ be a finite extension of
$\Qp$ containing a primitive $p$-th root of~$1$, of ramification index $e$ and
residual degree~$f$~; put $e_1=e/(p-1)$ and $q=p^f$.  Let $M=K(\!\root
p\of{K^\times})$ be the maximal kummerian extension of $K$ of exponent~$p$.
Let us first show, using the orthogonality relation (\S1), that the valuation
of the different of $M|K$ is
$$
v_M({\goth D}_{M|K})=(1+pe_1)pq^{e}-b_{(e+1)}-1,\leqno{(1)}
$$
where $b_{(e+1)}$ is the biggest break in the ramification filtration in
the lower numbering on $G=\Gal(M|K)$~; the lower breaks $b_{(i)}$ are computed
in prop.~{3}.

The orthogonality relation basically says that the filtration in the upper
numbering $(G^n)_{n\in[-1,+\infty[}$ on $G$ is given by $G^{-1}=G$, $G^0=G^1$,
$$
G^n=\bar U_{pe_1-n+1}^\perp\quad\qquad(n\in[1,pe_1+1]), \leqno{(2)}
$$
with the convention that $\bar U_0=K^\times\!/K^{\times p}$~; in
particular, $G^{pe_1+1}=\{\Id_M\}$.  Here the orthogonal is taken with respect
to the Kummer pairing $K^\times\!/K^{\times p}\times G\rightarrow{}_p\mu$.

It follows from $(2)$, the fact that the pairing is perfect (which implies
that $\dim_\Fp\bar U_m+\dim_\Fp\bar U_m^\perp=2+ef$ for every~$m$), and our
knowledge of $\Card\bar U_m$ (\citer\locdisc(prop.~{42})) that, for
$n\in[-1,+\infty[$,
$$
\hbox{\bf(}G^{n}:G^{n+1}\hbox{\bf)}=\cases{
1&if\/ $n>p{e_1}$,\cr
p&if\/ $n=p{e_1}$,\cr
1&if\/ $n<p{e_1}$ and\/ $p|n$,\cr
p&if\/ $n=-1$,\cr
q&otherwise.\cr}
$$
In the notation from \S1, this information can be summarised in one line~:
$$
\{1\}\subset_1 G^{pe_1}\subset_f\cdots\subset_f 
G^{pi+1}=G^{pi}\subset_f\cdots\subset_f G^1=G^0\subset_1G,
$$
where $i$ is any integer in $[1,e_1[$, and ``\thinspace$\subset_r$\thinspace''
means ``\thinspace codimension-$r$\thinspace''.

Thus, for $n\in[1,pe_1]$, we have $\dim_{\Fp}G^n
=1+\left(e-n+1+\left\lfloor{n-1\over p}\right\rfloor\right)\!f$~;
cf.~\citer\locdisc(prop.~{43}).  Here $e=(p-1)e_1$ is the ramification index,
and $f$ the residual degree, of $K|\Qp$.  In particular, $G^1$
(resp.~$G^{pe_1}$) has order $pq^e$ (resp.~$p$).

The upper ramification breaks occur therefore at $-1$, at the $e$ integers
$1=b^{(1)}<b^{(2)}<\cdots<b^{(e)}=pe_1-1$ in $[1,pe_1]$ which are prime
to~$p$, and at $b^{(e+1)}=pe_1$.  The order of the group drops by a factor of
$p$ at $-1$, by a factor of $q=p^f$ at each of the $b^{(i)}$ ($i\in[1,e]$),
and by a factor of $p$ at $pe_1$.  Consequently, we have the following table
for the index of $G^n$ in $G^0$~:
$$
\vbox{\halign{&\hfil$#$\hfil\quad\cr
n\in&\{0,1\}&]b^{(1)},b^{(2)}]&]b^{(2)},b^{(3)}]&
\cdots&]b^{(e-1)},b^{(e)}]&\{pe_1\}\cr
\noalign{\vskip-5pt}
\multispan7\hrulefill.\cr
\hbox{\bf(}G^{0}:G^{n}\hbox{\bf)}=&1&q&q^2&\cdots&q^{e-1}&q^e\cr
}}
$$
The lower ramification breaks occur therefore at $-1$ and at the $e+1$ integers
$b_{(i)}=\psi(b^{(i)})$ ($i\in[1,e+1]$), where $\psi=\psi_{M|K}$ is the
function on $[-1,+\infty[$ satisfying
$$
\psi(w)=\int_0^w\!\!\hbox{\bf(}G^0:G^t\hbox{\bf)}\,dt \leqno{(3)}
$$
\citer\serre(p.~74).  In view of the above table, it follows that $b_{(1)}=1$
and that $b_{(i+1)}=b_{(i)}+(b^{(i+1)}-b^{(i)})q^i$ for $i\in[1,e]$.  This may
also be verified using the formula $g_0\phi(r)=g_1+g_2+\cdots+g_r$
\citer\serre(p.~73), where $g_n=\Card G_n$.  

% Indeed, 
% $$\eqalign{
% g_0\phi(b_{i+1})&=g_0b^i+(b_{i+1}-b_i)g^{b^{i+1}}\cr
% &=g_0b^i+(b^{i+1}-b^i)p^{if}g^{b^{i+1}}\cr
% &=g_0b^i+(b^{i+1}-b^i)g_0,\cr
% }$$
% hence $\phi(b_{i+1})=b^{i+1}$.

(Notice that the $b_{(i)}$ are all $\equiv1\pmod p$, cf.~\citer\serre(p.~70)).

The $b_{(i)}$ can be computed recursively, starting from $b_{(1)}=1$.
Explicitly, for $i\in[1,e]$, we have
% $$
% b^{i+1}-b^i=\cases{
% 1&if\/ $i=e$ or if\/ $b^i\not\equiv-1\pmod p$\cr
% 2&if\/ $i\neq e$ and $b^i\equiv-1\pmod p$,\cr
% }$$
$$
b^{(i+1)}-b^{(i)}=\cases{
1&if\/ $i=e$ or if\/ $i\not\equiv0\mod{(p-1)}$\cr
2&if\/ $i\neq e\,$ and $\,i\equiv0\mod{(p-1)}$.\cr
}$$
(Notice that $(p-1)\,|\,i$ is equivalent to $p\,|\,(b^{(i)}+1)$.)
Therefore we have
$b_{(i)}=(1+q+\cdots+q^{i-1})+(q^{p-1}+\cdots+q^{a(i)(p-1)})$, with
$a(i)$ the integral part of $(i-1)/(p-1)$, for $i\in[1,e]$, and
$b_{(e+1)}=b_{(e)}+q^e$.  These are the expressions obtained in
\citer\delcorso(p.~287), albeit in the special case when
$K=F(\!\root{p-1}\of{F^\times})$ for some (finite) extension $F|\Q_p$.

To compute the valuation of the different ${\goth D}_{M|K}$ of $M|K$, it now
suffices to recall that the order $g_n$ of the ramification subgroup $G_n$ is 
$$
\vbox{\halign{&\hfil$#$\hfil\quad\cr
n\in&\{0,1\}&]b_{(1)},b_{(2)}]&]b_{(2)},b_{(3)}]&
\cdots&]b_{(e-1)},b_{(e)}]&]b_{(e)},b_{(e+1)}]\cr
\noalign{\vskip-5pt}
\multispan7\hrulefill.\cr
g_n=&pq^e&pq^{e-1}&pq^{e-2}&\cdots&pq&p\cr
}}
$$
As the valuation of the different ${\goth D}_{M|K}$ is
$\sum_{n\in\N}(g_n-1)$ \citer\serre(p.~64), we get
$$
v_M({\goth D}_{M|K})=(1+pe_1)pq^{e}-b_{(e+1)}-1,
$$
the expression claimed in $(1)$.  This expression follows from --- and indeed
led to --- the following lemma.
\th LEMMA {2}
\enonce
Let\/ $E|F$ be a finite galoisian extension of local fields, of group\/
$G=\Gal(E|F)$.  Suppose that the filtration\/ 
$(G^w)_{w\in[-1,+\infty[}$ has\/ $m>0$ positive breaks\/ $b^{(1)}, \ldots,
b^{(m)}$.  Let\/ $b_{(i)}=\psi_{E|F}(b^{(i)})$, for\/ $i\in[1,m]$, be the
breaks in the lower numbering.  Then the valuation of the different of\/ $E|F$
is\/ $v_E({\goth D}_{E|F})=(1+b^{(m)})g_0-(1+b_{(m)})$, where\/ $g_0=\Card
G_0$ is the order of the inertia subgroup. 
% $v_E({\goth  D}_{E|F})=(1+b_{(1)}+b^{(m)}-b^{(1)})g_0-b_{(m)}-1$, 
\endth 
For $w\in\,[0,+\infty[$, denote the index
$\hbox{\bf(}G^{0}:G^{w}\hbox{\bf)}$ of $G^w$ in $G^0$ as follows~:
$$
\vbox{\halign{&\hfil$#$\hfil\quad\cr
w\in&[0,b^{(1)}]&]b^{(1)},b^{(2)}]&\cdots&]b^{(m-1}),b^{(m)}]\cr
\noalign{\vskip-5pt}
\multispan5\hrulefill.\cr
\hbox{\bf(}G^{0}:G^{w}\hbox{\bf)}=&1&h^{(1)}&\cdots&h^{(m-1)}\cr
}}
$$
We have $b_{(i+1)}=b_{(i)}+(b^{(i+1)}-b^{(i)})h^{(i)}$ for every
$i\in\,[1,m[$.  The cardinality of the subgroups $G_t$ ($t\in[0,+\infty[$) of
$G_0$ in the lower numbering is~:
$$
\vbox{\halign{&\hfil$#$\hfil\quad\cr
t\in&[0,b_{(1)}]&]b_{(1)},b_{(2)}]&\cdots&]b_{(m-1)},b_{(m)}]\cr
\noalign{\vskip-5pt}
\multispan5\hrulefill.\cr
\Card G_{t}=&g_0&g_0/h^{(1)}&\cdots&g_0/h^{(m-1)}\cr
}}
$$
Now, for $i\in[1,m[$, the contribution of the interval $]b_{(i)},b_{(i+1)}]$ to
the sum $\sum_{n\in\N}(\Card G_n-1)$ is
$$\eqalign{
(b_{(i+1)}-b_{(i)})\left({g_0-h^{(i)}\over h^{(i)}}\right)
&=(b^{(i+1)}-b^{(i)})h^{(i)}\left({g_0-h^{(i)}\over h^{(i)}}\right)\cr
&=(b^{(i+1)}-b^{(i)})\left(g_0-h^{(i)}\right),\cr
}$$
and hence $v_E({\goth D}_{E|F})$, the sum over the contributions of these
$m-1$ intervals $]b_{(i)},b_{(i+1)}]$ ($i\in[1,m[$) and of the initial segment
$[0,b_{(1)}]$ of $1+b_{(1)}$ points is given by  
% $$
% \eqalign{
% &(1+b_{(1)})(g_0-1)+(b^{(2)}-b^{(1)})(g_0-h^{(1)})+
%  \cdots+(b^{(m)}-b^{(m-1)})(g_0-h^{(m-1)})\cr
% &=
% }$$
$$
(1+b_{(1)}+b^{(m)}-b^{(1)})g_0-b_{(m)}-1
=(1+b^{(m)})g_0-(1+b_{(m)}),
$$
proving lemma~{2}.  (The case $m=1$ implies \citer\locdisc(prop.~{60}), in
view of fact that for an $\Fp$-line $D\subset\bar U_c$, $D\not\subset\bar
U_{c+1}$, the unique ramification break of the kummerian extension $K(\!\root
p\of D)|K$ occurs at $b^{(1)}=b_{(1)}=pe_1-c$ if $c\neq pe_1$.)

{\it Example}\pointir Take $F=\Qp$ and $E=\Qp(\xi_a)$, where $\xi_a$ is a
primitive $p^a$-th root of~$1$ for some $a>0$.  When $p\neq2$, we have
$m=a$, $b^{(1)}=b_{(1)}=0$, $b^{(m)}=a-1$, $b_{(m)}=p^{a-1}-1$ and
$g_0=\varphi(p^a)$.  Therefore $v_E({\goth
  D}_{E|F})=a\varphi(p^a)-p^{a-1}$.  Consider now $p=2$.  If $a=1$, the
extension $E|F$ is trivial and the lemma does not apply~; nor do we need
to apply it.  For $a>1$, the only change is that $m=a-1$,
$b^{(1)}=b_{(1)}=1$, leading to the same result~: $v_E({\goth
  D}_{E|F})=a\varphi(2^a)-2^{a-1}=(a-1)2^{a-1}$.

Let us summarise what we have learnt about our maximal kummerian extension
$M|K$ of exponent~$p$.  Let $\displaystyle a(i)=\left\lfloor{i-1\over
    p-1}\right\rfloor$.

\th PROPOSITION {3} 
\enonce
The\/ $e+1$ positive ramification breaks of\/ $M|K$ occur at\/ $b^{(i)}=i+a(i)$
for $i\in[1,e]$, and at\/ $b^{(e+1)}=pe_1$, in the upper numbering.  In
the lower numbering, they occur at
$$
b_{(i)}=(1+q+\cdots+q^{i-1})+(q^{p-1}+\cdots+q^{a(i)(p-1)}) \qquad
(i\in[1,e]) \leqno{(4)}
$$
and at\/ $b_{(e+1)}=b_{(e)}+q^e$.  We have\/ $v_M({\goth
  D}_{M|K})=(1+pe_1)pq^e-(1+b_{(e+1)})$ and\/ $v_K(d_{M|K})=p.v_M({\goth
  D}_{M|K})$. 
\endth 
The statement about the discriminant follows from the fact that the residual
degree of $M|K$ is $p$.  Notice that $a(e)=e_1-1$, so that
$b_{(e+1)}=
(1+q+q^2+\cdots+q^{e})+(q^{p-1}+q^{2(p-1)}+\cdots+q^{(e_1-1)(p-1)})$, with
$q=p^f$.

{\it Example}\pointir Take $K=\Qp(\zeta)$, where $\zeta$ is a primitive $p$-th
root of~$1$. Then $e=p-1$, $e_1=1$, and $f=1$.  The $p$ ramification breaks of
$M|K$ are $1$, $2$, $\ldots$, $p$ in the upper numbering~; $1$, $1+p$,
$\ldots$, $1+p+p^2+\cdots+p^{p-1}$ in the lower numbering.  Therefore
$v_M({\goth D}_{M|K})=p^{p+1}+p^p-p^{p-1}-\cdots-p-2$.

{\it Example}\pointir The last result of \citer\delcorso(p.~287) can be
recovered by taking $K=F(\!\root p-1\of{F^\times})$, where $F$ is any finite
extension of $\Qp$.  Keep the notation $e=(p-1)e_1$ and $q=p^f$ relative to
$K$.  We have $v_K({\goth D}_{K|F})=p-2$, so the valuation of ${\goth
  D}_{M|F}={\goth D}_{M|K}{\goth D}_{K|F}$ \citer\serre(p.~51) is
$$\eqalign{
v_M({\goth
  D}_{M|F})&=v_M({\goth D}_{M|K})+v_M({\goth D}_{K|F})\cr
&=((1+pe_1)pq^e-b_{(e+1)}-1)+pq^e(p-2)\cr
&=(e_1+1)p^2q^e-pq^e-(1+b_{(e+1)}).\cr
%&=(e_1+1)p^2q^e-pq^e-q^e-{q^e-1\over q-1}-{q^e-1\over q^{p-1}-1}.\cr
}$$

It is also possible to deduce the following result of J.-M. Fontaine.
\th COROLLARY {4} (\citer\fontaine(p.~362))
\enonce
Let\/ $F|\Qp$ be a finite extension, $E|F$ a totally ramified elementary
abelian\/ $p$-extension.  Every upper ramification break\/ $u$ of\/ $E|F$ is
in\/ $[1,pe_1]$ and is prime to\/~$p$, unless\/ $u=pe_1$.  The order of\/
$\Gal(E|F)^{pe_1}$ is\/ $1$ or\/ $p$.
\endth
One may add that if\/ $pe_1$ occurs for some\/ $E|F$, then\/ $F$ contains a
primitive\/ $p$-th root\/ $\zeta$ of\/~$1$ (cf.~\citer\locdisc(prop.~{63})) and
there is a uniformiser $\pi$ of $F$ such that $\pi\in E^{\times p}$.

One may ask for a converse~: which sequences do occur as the upper
ramification breaks of an elementary abelian\/ $p$-extension $E|F$~?  We may
ask for the maximal degree $[E:F]$ when there is a single break.  We may ask
for the number of extensions with given ramification breaks.  If $\zeta\in F$,
these questions can be answered by Kummer theory~; see below.  If
$\zeta\not\in F$, we may reduce to the previous case by considering the
extension $E(\zeta)$ of $F(\zeta)$, as in the proof of
\citer\locdisc(prop.~{63}).  Alternatively, one may appeal to local class
field theory, to which we turn in a moment.

{\it The existence of exponent-$p$ kummerian\/ extensions with given upper
  ramification breaks.}  Suppose that $K$ is a finite extension of
$\Qp(\zeta)$.  Every strictly increasing sequence $u_1<u_2<\ldots<u_n$ ($n>0$)
of numbers which are in $[1,pe_1]$, with the possible exception of $u_1$,
which can be $-1$, and which are all prime to~$p$, with the possible exception
of $u_n$, which can be $pe_1$, is the sequence of upper ramification breaks of
some exponent-$p$ kummerian\/ extension $L$ of $K$.  

Note first that we need only consider the case $u_1\neq-1$.  For if $L_1$ is
the unique unramified degree-$p$ extension of $K$, and if $L_2$ is a (totally
ramified) exponent-$p$ kummerian extension with upper ramification breaks
$u_2,\ldots,u_n$, then $L=L_1L_2$ has ramification breaks $-1,u_2,\ldots,u_n$
in the upper numbering.  Assume therefore that $u_1>0$.

We may look for $L$ inside the maximal exponent-$p$ kummerian extension
$M=K(\!\root p\of{K^\times})$ of $K$.  Equivalently, we look for a
subgroup $H$ of $G=\Gal(M|K)$ such that
$$
\setbox0\hbox{\phantom{$\neq(G/H)^{u_2+1}=\cdots=(G/H)^{u_3}$}}
\setbox1\hbox to\wd0{\dotfill}
\eqalign{G/H=(G/H)^{-1}=\cdots=(G/H)^{u_1}
&\neq(G/H)^{u_1+1}=\cdots=(G/H)^{u_2}\cr
&\neq(G/H)^{u_2+1}=\cdots=(G/H)^{u_3}\cr
&\box1\cr
&\neq(G/H)^{u_{n-1}+1}=\cdots=(G/H)^{u_n}\cr
&\neq(G/H)^{u_n+1}=\{\bar1\}\cr
}$$
and take $L=M^H$.  In view of the compatibility
of the upper-numbering filtration with the passage to the quotient, we have
$(G/H)^i=G^iH/H$, and we are led to require 
$$
\setbox0\hbox{\phantom{$\neq G^{u_2+1}H=\cdots=G^{u_3}H$\quad}}
\setbox1\hbox to\wd0{\dotfill}
\eqalign{G=G^{-1}H=\cdots=G^{u_1}H
&\neq G^{u_1+1}H=\cdots=G^{u_2}H\cr
&\neq G^{u_2+1}H=\cdots=G^{u_3}H\cr
&\box1\cr
&\neq G^{u_{n-1}+1}H=\cdots=G^{u_n}H\cr
&\neq G^{u_n+1}H=H.\cr
}$$
If we identify $G$ with $\bar U_0=K^\times\!/K^{\times p}$ using the
reciprocity isomorphism, and recall the structure of the filtered $\Fp$-space
$\bar U_0$ (\citer\locdisc(prop.~{42})), we may conclude that such a subspace
$H\subset G$ exists always.

However, this appeal to local class field theory can be avoided when $\zeta\in
K$, as here.  Appeal can be made instead to the orthogonality relation $(2)$,
$G^{n\perp}=\bar U_{pe_1-n+1}$ for $n\in[1, pe_1+1]$.  So we look for a
subspace $D\subset\bar U_0$ such that, writing $D_i$ for $\bar U_i\cap D$, we
have
$$
\setbox0\hbox{\phantom{$\neq D_{pe_1-u_1}=\cdots=D_{pe_1-u_2+1}$\ \ }}
\setbox1\hbox to\wd0{\dotfill}
\eqalign{\{\bar1\}=D_{pe_1}=\cdots=D_{pe_1-u_1+1}
&\neq D_{pe_1-u_1}=\cdots=D_{pe_1-u_2+1}\cr
&\neq D_{pe_1-u_2}=\cdots=D_{pe_1-u_3+1}\cr
&\box1\cr
&\neq D_{pe_1-u_{n-1}}=\cdots=D_{pe_1-u_n+1}\cr
&\neq D_{pe_1-u_n}=\cdots=D_0=D\cr
}$$
and take $H=D^\perp$.  In view of \citer\locdisc(prop.~{42}), it is clear that
such a $D$ exists always.  It is also clear how to get every such $D$, and how
to count the number of such~$D$.  (Notice that every such $D$ is an
$\Fp$-point of a certain open subvariety of the $\Fp$-variety of all subspaces
of $\bar U_0$.)

More explicitly, choose, for every $i\in[1,n]$, a $d_i$-dimensional ($d_i>0$)
sub-$\Fp$-space $E_i\subset\bar U_{pe_1-u_i}$ whose intersection with $\bar
U_{pe_1-u_i+1}$ is $\{\bar1\}$~; such choices are possible in view of
\citer\locdisc(prop.~{42}).  Finally take $D=E_1E_2\cdots E_n$, so that
$D_{pe_1-u_i}=E_1E_2\cdots E_i$.  The dimension of $E_i$ is between $1$ and
$f$, except when $i=n$ and $u_n=pe_1$, where $\dim_{\Fp} E_{pe_1}=1$.

% {\it The existence of kummerian\/ exponent-$p$ extensions with a given
%   sequence of upper ramification breaks.}  Suppose that $K$ is a finite
% extension of $\Qp(\zeta)$.  Every strictly increasing sequence
% $u_1<u_2<\ldots<u_n$ of numbers in $[1,pe_1]$ which are all prime to~$p$,
% except possibly $u_n$, which can be $pe_1$, is the sequence of upper
% ramification breaks of some totally ramified kummerian\/ exponent-$p$
% extension $L$ of $K$.  To see this, choose, for every $i\in[1,n]$, a
% sub-$\Fp$-space $D_i\subset\bar U_{pe_1-u_i}$ whose intersection with $\bar
% U_{pe_1-u_i+1}$ is $\{\bar1\}$~; such choices are possible in view of
% prop.~{42}.  Then the upper ramification breaks of the extension $L=K(\root
% p\of{D_1D_2\cdots D_n})$ are precisely $u_1,u_2,\ldots,u_n$.
 
{\it The degree of a totally ramified exponent-$p$ kummerian\/ extension with
  a single ramification break.}  Let $L|K$ be such an extension, and let $u$
be the unique ramification break.  As we have seen, $[L:K]=p^j$, where $j=1$
if $u=pe_1$ and $j\in[1,f]$ if $u\neq pe_1$~; in the latter case, there are
extensions with any preassigned $j\in[1,f]$.  The valuation of the
discriminant is $v_K(d_{L|K})=(p^j-1).(1+u)$~; the case $j=1$ is
\citer\locdisc(cor.~{64}).
  
{\it The number of exponent-$p$ kummerian\/ extensions with a given sequence
  of upper ramification breaks.}  It is clear that every such extension $L|K$
with given upper ramification breaks $u_1<u_2<\ldots<u_n$ as above arises as
$L=K(\!\root p\of{E_1E_2\cdots E_n})$ for some choice of subspaces $E_i$.  For a
different choice $E_i^\prime$ of the subspaces, we get the same extension if
and only if $E_1E_2\cdots E_n=E_1^\prime E_2^\prime\cdots E_n^\prime$.  This
leads to a mildly complicated counting problem --- how many subspaces
$D\subset\bar U_0$ are there such that $D_i=D\cap\bar U_i$ satisfy the
conditions displayed above~?  --- which can be solved in any given instance~;
see \citer\locdisc(cor.~{66}) for the degree-$p$ cyclic case.  At the other
extreme, when $n=e+1$ and $d_i$ are as large as they can be, namely $d_i=f$
for $i\in[1,e]$ and $d_{e+1}=1$, a subspace $D\subset\bar U_0$ is a solution
if and only if it is a hyperplane not containing the line $\bar U_{pe_1}$.
  
{\it Remark}\pointir Let $D,D'$ be two distinct lines in $\bar U_b$ neither of
which is contained in $\bar U_{b+1}$.  If the plane $DD'$ meets $\bar U_{b+1}$
only at the origin $\{\bar1\}$, then the extension $L=K(\!\root p\of{DD'})$ has
a single ramification break, namely $pe_1-b$.  Otherwise, let $a$ be the
largest integer such that $DD'\cap\bar U_{a}\neq\{\bar1\}$~; we have
$a\in\,]b,pe_1]$.  The extension $L|K$ now has two ramification breaks, namely
$pe_1-a$ and $pe_1-b$ if $a\neq pe_1$, and $-1$ and $pe_1-b$ if $a=pe_1$.  For
an extreme example, where $b=0$ and $a=pe_1$, see \citer\locdisc(ex.~{51}).

{\it The valuation of the different of an exponent-$p$ kummerian extension
  with given upper ramification breaks.}  Let $u_1<u_2<\ldots<u_n$ be a
strictly increasing sequence of numbers in $[1,pe_1]$ which are all prime
to~$p$ except possibly $u_n$, which can be $pe_1$~; choose $E_i$ as above and
let $L=K(\!\root p\of{E_1E_2\cdots E_n})$~; the upper ramification breaks of
$L|K$ occur at $u_1,u_2,\ldots,u_n$.  Let $d_i$ be the dimension of $E_i$, so
that $\Card E_i=p^{d_i}$.  The valuation of the different of the extension
$L|K$ is given by $v_L({\goth
  D}_{L|K})=(1+u_n)p^{d_1+d_2+\cdots+d_n}-(1+\psi_{L|K}(u_n))$ (lemma~{2}),
where
$$
\psi_{L|K}(u_n)=u_1.1+
(u_2-u_1).p^{d_1}+
\cdots+
(u_n-u_{n-1}).p^{d_1+d_2+\cdots d_{n-1}},
$$
as follows from the definition of $\psi_{L|K}$ $(3)$ and the following piece of
information about the ramification filtration on $G=\Gal(L|K)$~:
$$
\vbox{\halign{&\hfil$#$\hfil\quad\cr
j\in&[0,u_1]&]u_1,u_2]&]u_2,u_3]&\cdots&]u_{n-1},u_n]\cr
\noalign{\vskip-5pt}
\multispan6\hrulefill.\cr
\hbox{\bf(}G^{0}:G^{j}\hbox{\bf)}
=&1&p^{d_1}&p^{d_1+d_2}&\cdots&p^{d_1+d_2+\cdots d_{n-1}}\cr 
}}
$$
Notice that $v_L({\goth D}_{L|K})$ depends on the subspaces $E_1$, $E_2$,
$\ldots$, $E_n$ only via the breaks $u_i$ and the dimensions $d_i$.

{\it Remark}\pointir This can be used to compute the minimum or the maximum
value of $v_K(d_{L|K})$, where $L$ runs through totally ramified elementary
abelian $p$-extensions of $K$ of given degree $p^m$ ($m>0$).  Suppose first
that we want to maximise $v_K(d_{L|K})$.  If $m=1$, we would take $L=K(\root
p\of D_1)$, where $D_1$ is a line in $\bar U_0$ which is not in $\bar U_1$.  If
$m=2$, we would take $L=K(\root p\of{D_1D_2})$, where $D_2$ is a line in $\bar
U_1$ not in $\bar U_2$. The idea is the same for higher $m$.

When $p=2$ and $K|\Q_2$ is unramified, Tate (letter to Serre, July 1973) gives
the upper bound $3.2^{m-1}+2^m-2$ for $v_K(d_{L|K})$ using the formula
expressing the discriminant as the product of the conductors of characters of
$\Gal(L|K)$ (the {\it F{\"u}hrerdiskriminantenproduktformel\/}) and local
class field theory to compute the conductors \citer\tate(p.~155).

Suppose next that we want to minimise $v_K(d_{L|K})$.  If $m=1$, we would take
$L=K(\root p\of D_1)$, where $D_1\neq\bar U_{pe_1}$ is any line in $\bar
U_{pe_1-1}$.  If $m=2$ and $f>1$, we would take $L=K(\root p\of{D_1D_2})$,
where $D_2\neq\bar U_{pe_1}$ is any line in $\bar U_{pe_1-1}$ distinct from
$D_1$.  It is easy to guess what to do for $m=2$ and $f=1$, and ideed for any
$m,f$.

{\it Example}\pointir For an example of the smallest possible degree having
every possible upper ramification break, take $n=e+1$, $u_i=b^{(i)}=i+a(i)$
for $i\in[1,e]$, $u_{e+1}=pe_1$, and $d_i=1$ for every $i$.  In the lower
numbering, the breaks occur at
$$
l_i=(1+p+\cdots+p^{i-1})+(p^{p-1}+\cdots+p^{a(i)(p-1)}) \qquad
(i\in[1,e]) 
$$
and at\/ $l_{e+1}=l_e+p^e$.  We have\/ $v_L({\goth
  D}_{L|K})=(1+pe_1)p^{e+1}-(1+l_{e+1})$ (lemma~{2}).

Let us summarise a part of our discussion of exponent-$p$ kummerian extensions
with given ramification breaks.

\th PROPOSITION {5} 
\enonce
Every strictly increasing sequence\/ $u_1<u_2<\ldots<u_n$ of integers in\/
$[1,pe_1]$ which are all prime to\/~$p$, with the possible exception of\/
$u_n$, which can be\/ $pe_1$, is the sequence of upper ramification breaks of
some (totally ramified) exponent-$p$ kummerian\/ extension\/ $L$ of\/ $K$.  If
a break occurs at\/ $pe_1$, then there is a uniformiser\/ $\pi$ of\/ $K$ such
that\/ $\pi\in L^{\times p}$ and conversely.  When there is a single break\/
$u$, there are $f$ possibilities for the degree\/ $[L:K]$ if\/ $u\neq pe_1$,
namely $p,p^2,\ldots,p^f$, but only one possibility, namely\/ $[L:K]=p$, if\/
$u=pe_1$. 
\endth
As noted earlier, the ultimate source of this prop.\ is the analysis of the
filtered $\Fp$-space $K^\times\!/K^{\times p}$ (as for example in
\citer\locdisc()), whether we use local class field theory or Kummer theory
(as here).

\medskip \centerline{***} \medbreak

Now suppose that $F|\Qp$ is a finite extension for which $F^\times$ does not
have an element of order~$p$ (and hence $p\neq2$), let $N$ be the maximal
abelian extension of~$F$ of exponent~$p$, and let $G=\Gal(N|F)$.  Local class
field theory provides an isomorphism $F^\times\!/F^{\times p}\rightarrow G$
under which $\bar U_n$ surjects onto $G^n$ for every $n>0$~; thus the inertia
group $G^0=G^1$ has order $q^e$ and index $p$ in $G$.  There are now only $e$
positive ramification breaks \citer\locdisc(prop.~42)~; in the upper
numbering, they occur at the $e$ integers $b^{(i)}=i+a(i)$ ($i\in[1,e]$) in
the interval $[1,pe_1[$ which are prime to~$p$.  At each break, the order of
the group drops by a factor of $q$, so the $e$ lower ramification breaks are
as given in $(4)$.  The valuation of the different is $v_N({\goth
  D}_{N|F})=(1+b^{(e)})q^e-(1+b_{(e)})$ (lemma~{2}).

Here and elsewhere, we have made use of the following elementary fact~: for
every prime~$p$ and every integer~$m>0$, the~$c(m)=m-\lfloor m/p\rfloor$
integers in $[1,m]$ which are prime to~$p$ constitute the image of the
strictly increasing function $b^{(\phantom{i})}:[1,c(m)]\rightarrow[1,m]$,
$b^{(i)}=i+a(i)$.
 
{\it Example}\pointir Take $F=\Qp$ ($p\neq2$).  There is a unique positive
ramification break, at $b^{(1)}=b_{(1)}=1$ (see the parenthetical remark
before cor.~{64} of \citer\locdisc()), so $v_N({\goth D}_{N|F})=2(p-1)$.
(When $p=2$, $N$ is the maximal kummerian extension of $\Q_2$, which has been
treated earlier.)

Let $K$ be a finite extension of $\Qp$.  In principle, the determination of
the functions $\phi_{L|K}$, $\psi_{L|K}$, for any finite galoisian extension
$L|K$ can be reduced to the case treated in \citer\locdisc(prop.~{63}) (cyclic
of degree~$p$), and indeed to the case treated in \citer\locdisc(cor.~{62})
(kummerian of degree~$p$).  But a little bit of local class field theory makes
life simpler.

{\it Example}\pointir Let $K$ be any finite extension of $\Qp$, let $m>0$ be
an integer, let $\pi$ be a uniformiser of $K$, let $L$ be the unique totally
ramified abelian extension of~$K$ such that $\pi\in N_{L|K}(L^\times)$ and
$N_{L|K}(\ogoth_L^\times)=U_m$, and let $G=\Gal(L|K)=\ogoth_K^\times/U_m$.
(If $K=\Qp$ and $\pi=p$, then $L=\Qp(\xi_m)$, which we have treated in~\S~1.)
Then $G^{-1}=G^0$ has order $q^{m-1}(q-1)$, and the index of $G^n$ in $G^0$
is~:
$$
\vbox{\halign{&\hfil$#$\hfil\quad\cr
n=&0&1&2&\cdots&m\cr
\noalign{\vskip-5pt}
\multispan6\hrulefill.\cr
\hbox{\bf(}G^{0}:G^{n}\hbox{\bf)}=&1&q-1&q(q-1)&\cdots&q^{m-1}(q-1)\cr
}}
$$
If the residual cardinality $q$ of~$K$ is $\neq2$, there are $m$
ramification breaks, $b^{(i)}=i-1$ for $i\in[1,m]$, in the upper numbering~;
in the lower numbering, they are $b_{(1)}=0$, $b_{(2)}=q-1$, $\ldots$,
$b_{(m)}=q^{m-1}-1$. In particular, $v_K(d_{L|K})=mq^{m-1}(q-1)-q^{m-1}$,
which is independent of $\pi$.

If $q=2$ and $m>1$, there are only $m-1$ breaks, at $b^{(i)}=i$ for
$i\in[1,m[$.  If $q=2$ and $m=1$, the extension is trivial.  The expression
for the valuation of the discriminant remains the same in all cases~;
cf.~\citer\iwasawa(p.~110).
 
{\it Remark}\pointir Take $m=1$, so that $L|K$ is cyclic of degree~$q-1$, and
hence, being totally ramified, obtained by adjoining $\root{q-1}\of\varpi$ for
some uniformiser~$\varpi$, uniquely determined up to $1$-units.  Which
uniformiser~?  The answer for $K=\Qp$ and $\pi=p$ is $\varpi=-p$
\citer\locdisc(prop.~{24}).  We have $\varpi=-\pi$ in general, for $L$ is the
splitting field of $T^q+\pi T$ \citer\neukirch(p.~61).  Turning things around,
we may say that if $L=K(\!\root{q-1}\of\varpi)$ for some uniformiser $\varpi$
of $K$, then $N_{L|K}(L^\times)/K^{\times q-1}$ is the subgroup of
$K^\times\!/K^{\times q-1}$ generated by the image of $-\varpi$ (which is the
same as the image of $\varpi$ when $q$ is even, for then $-1=(-1)^{q-1}$).  It
reflects the fact that the product of all elements in the multiplicative group
$k^\times$ of the residue field is $-1$~; indeed, the conjugates of
$\root{q-1}\of\varpi$ are $u.\root{q-1}\of\varpi$, as $u$ runs through
$k^\times$.

Recall that the Hilbert symbol $K^\times\!/K^{\times q-1}\times
K^\times\!/K^{\times q-1}\rightarrow{}_{q-1}\mu$ has the property that if
$L=K(\!\root{q-1}\of D)$ for some subgroup $D\subset K^\times\!/K^{\times
  q-1}$, then $D^\perp=N_{L|K}(L^\times)/K^{\times
  q-1}$~\citer\fesvost(p.~144).  This means that
$D_{\varpi}^\perp=D_{-\varpi}$, where $D_a\subset K^\times\!/K^{\times q-1}$
is the subgroup generated by $a\in K^\times$.  It might be added that for
$D=\ogoth^\times\!/\ogoth^{\times q-1}=k^\times$, we have $D^\perp=D$, because
$K(\!\root {q-1}\of D)$ is the unramified extension of~$K$ of degree~$q-1$.

{\it An application of the\/} F{\"u}hrerdiskriminantenproduktformel\pointir
Let us show how the discriminant of the maximal kummerian extension $M|K$ of
exponent~$p$ could have been computed by an application of this formula after
we had determined the possible ramification breaks $t$ of a degree-$p$
kummerian extension $L|K$ \citer\locdisc(cor.~{62}) and the number of
extensions for which a given break occurs \citer\locdisc(cor.~{66}), if we knew
that the exponent of the conductor of $L|K$ is $t+1$.  Class field theory
provides this last bit of knowledge, for it says that, under the reciprocity
map, $U_t$ surjects onto $\Gal(L|K)$ whereas the image of $U_{t+1}$ is
$\{\Id_M\}$, so $t+1$ is the smallest integer $m$ such that $U_m\subset
N_{L|K}(L^\times)$.

Recall that the formula in question, applied to an abelian extension $E|F$ of
local fields, says that the discriminant ideal ${\goth d}_{E|F}$ equals
$\prod_\chi{\goth f}(\chi)$, where the product is taken over all characters
$\chi:\Gal(E|F)\rightarrow\C^\times$ and ${\goth f}(\chi)$ is the
conductor of $\chi$ \citer\iwasawa(p.~113), \citer\serre(p.~104)~; of course,
only ramified characters need be considered.

To a ramified character $\chi$ of $G=\Gal(M|K)$ corresponds a ramified
degree-$p$ cyclic extension $L=M^{\Ker(\chi)}$, and each ramified degree-$p$
cyclic extension $L$ arises from $p-1=\Card\Aut({}_p\mu)$ characters~$\chi$.
In view of this, it is sufficient to compute $(p-1)\sum_{i=1}^{e+1}
(t_i+1).n_i$, where
$$
t_i=i+a(i) \quad(i\in[1,e])~;\quad t_{e+1}=pe_1
$$
are the possible positive ramification breaks \citer\locdisc(cor.~{62})
and, as shown in \citer\locdisc(cor.~{66}),
$$
n_i=p^{(i-1)f+1}+\cdots+p^{if}={p(q^i-q^{i-1})\over p-1}\quad(i\in[1,e])~;
\quad n_{e+1}=pq^e\qquad 
$$
with $q=p^f$, is the number of ramified degree-$p$ cyclic extensions of $K$
whose ramification break occurs at~$t_i$.  Now, it is easily seen that
$$\eqalign{
&(p-1)\sum_{i=1}^e n_i=p.\left(q^e-1\right),\quad
(p-1)\sum_{i=1}^e in_i=p.\left(eq^e-{q^e-1\over q-1}\right),\cr
&\qquad\quad(p-1)\sum_{i=1}^e a(i) n_i
  =p.\left((e_1-1)q^e-{q^e-1\over q^{p-1}-1}+1\right),\cr
}$$
where, to compute the last sum,  rewrite it as
$\sum_{i=1}^e=\sum_{j=0}^{e_1-1}\sum_{r=1}^{p-1}$ and recall that $a(i)=j$ if
$i=(p-1)j+r$ for some $j\in\N$ and $r\in[1,p-1]$.
% recall that $a(i)$ is the integral part of
% $(i-1)/(p-1)$, and hence
% $$
% a(i+1)-a(i)=\cases{
% 0&if\/ $i\not\equiv0\mod{(p-1)}$\cr
% 1&if\/ $i\equiv0\mod{(p-1)}$.\cr
% }$$
Therefore $\displaystyle (p-1)\sum_{i=1}^{e}(t_i+1)n_i
=p.\left(e_1pq^e-{q^e-1\over q-1}-{q^e-1\over
    q^{p-1}-1}\right)$ and, by the 
{\it F{\"u}hrerdiskriminantenproduktformel}, $v_K(d_{M|K})$ equals
$$
(p-1)\sum_{i=1}^{e+1}(t_i+1)n_i
=p.\left((e_1p^2+p-1)q^e-{q^e-1\over q-1}-{q^e-1\over q^{p-1}-1}\right),
\leqno{(5)}
$$
which is the same as in prop.~{3}.  This computation can be encapsulated in
the following lemma.

\th LEMMA {6}
\enonce
Let\/ $p>1$, $e_1>0$ be integers and $q>1$ real~; put\/ $e=e_1.(p-1)$.  For\/ 
$i\in[1,e]$, let\/ $a(i)$ be the integral part of\/ $(i-1)/(p-1)$~; define\/
$t_i=i+a(i)$, $n_i=p(q^i-q^{i-1})/(p-1)$ for\/ $i\in[1,e]$, and\/
$t_{e+1}=pe_1$, $n_{e+1}=pq^e$.  Then\/ $(p-1)\sum_{i=1}^{e+1}(t_i+1).n_i$ is
given by\/~$(5)$.   
\endth
Maximal elementary abelian $p$-extensions can be treated in like manner.

It is not surprising that the {\it F{\"u}hrerdiskriminantenproduktformel\/} can
compute the discriminant without reference to the lower ramification
filtration.  Indeed, information about this filtration goes into the proof of
the formula \citer\iwasawa(p.~113).

{\it Totally ramified finite abelian $p$-extensions and the endomorphism of
  raising to the exponent~$p$.}  Let $F$ be a finite extension of $\Qp$ and
denote by $U_n$ ($n>0$) the groups of higher principal units of $F$.  As
always, $e$ is the absolute ramification index and $e_1=e/(p-1)$.

Let $G$ be a finite commutative $p$-group.  A result of Fontaine
\citer\fontaine(p.~362) about totally ramified $G$-extensions $E|F$ follows
from the study of the raising-to-the-exponent-$p$ map $(\ )^p$ on the
$\Zp$-modules $U_n$ and the fact that the reciprocity map $F^\times\rightarrow
G$ carries $U_n$ onto $G^n$.  As $E|F$ is totally ramified, $G=G^0$, and as
$G$ is a $p$-group, $G^0=G^1$, so $G$ is essentially a quotient of $U_1$ by a
(closed) subgroup of finite index.  

Recall that the map $(\ )^p$ carries $U_n$ into $U_{\lambda(n)}$
\citer\locdisc(prop.~{27}), where $\lambda(n)=\inf(pn, n+e)$, and that
$U_n^pU_{\lambda(n)+1}=U_{\lambda(n)}$ in all cases except when $F^\times$ has
an element of order~$p$ and $n=e_1$, in which case $U_{e_1}^pU_{pe_1+1}$ has
index~$p$ in $U_{pe_1}$ \citer\locdisc(prop.~{29}).  Hence the following
result, which sharpens \citer\fontaine(p.~362).  For a subgroup $A$ of $G$,
denote by $A^{(p)}$ the image of $A$ under the endomorphism $(\ )^p$, because
the notation $A^p$ is in conflict with the upper numbering.  (The result
$(G^u)^{(p)}\subset G^{\lambda(u)}$ continues to hold even if we allow the
residue field to be merely perfect \citer\sen(p.~45).)

\th PROPOSITION {7} 
\enonce
Let\/ $E|F$ be a totally ramified finite abelian\/ $p$-extension and let\/
$G=\Gal(E|F)$.  Then equality holds in\/ $(G^n)^{(p)}G^{\lambda(n)+1}\subset
G^{\lambda(n)}$ in all cases except possibly when\/ $n=e_1$~; for\/ $n=e_1$,
the index can be\/ $1$ or\/ $p$.  The index-$p$ case occurs precisely when\/
$H^{pe_1}\neq\{1\}$, where $H=G/G^{(p)}$ is the maximal elementary abelian
quotient of\/ $G$.
\endth

Write $L=E^{G^{(p)}}$, so that $H=\Gal(L|F)$.  As $H$ is an elementary abelian
$p$-group, we always have $H^{pe_1+1}=\{1\}$ (cor.~{4}).
Also, if $H^{pe_1}\neq\{1\}$, then $\zeta\in F$ and there is a uniformiser
$\pi$ of $F$ such that $\pi\in L^{\times p}$, and conversely (prop.~{5}).

Remark finally that, going modulo ${G^{(p)}}$, we may assume that $G$ is
elementary abelian, in which case $H=G$, and we are reduced to cor.~{4}, and
ultimately to \citer\locdisc(prop.~{42}), whether we use the orthogonality
relation $(2)$ or the reciprocity isomorphism.

\bigbreak

{\bf 5.  Summary of local Artin-Schreier theory.}  Let us first summarise our
results in the Artin-Schreier theory for local function fields of
characteristic~$p$.  These results were arrived at by analogy with the Kummer
theory of local number fields as recalled in \S1, and they may be considered
as a refinement of the theory presented in standard textbooks such as
\citer\fesvost(Chapter~III, \S2, p.~74).  The actual writing of \S6, which
contains the proofs, was spurred by a fortuitous encounter with
\citer\wuscheidler() and was achieved within a few days.

Let $k|\F_p$ be a finite extension, $f=[k:\Fp]$ its degree, and let $K$ be a
local field with $k$ as the field of constants (and the residue field).
Denote by $\ogoth$ the ring of integers of $K$, and by $\pgoth\subset\ogoth$
the unique maximal ideal of $\ogoth$.  If we choose a uniformiser $\pi$ of $K$
(which we don't need to), we have $K=k(\!(\pi)\!)$, $\ogoth=k[[\pi]]$,
$\pgoth=\pi\ogoth$.

The filtration $(\pgoth^n)_{n\in\Z}$ on the additive group $K$ by powers of
$\pgoth$ induces a filtration on the $\Fp$-space $\overline{K}=K/\wp(K)$,
where $\wp$ is the endomorphism $x\mapsto x^p-x$ of $K$.  We denote the
induced filtration by $(\overline{\pgoth^n})_{n\in\Z}$~; we have
$\overline{\pgoth}=\{0\}$ (lemma~{8}), and the codimension at each step is
given by 
$$
\{0\}
\subset_1\overline\ogoth
\subset_f\overline{\pgoth^{-1}}
\subset_f\cdots
\subset_f\overline{\pgoth^{-pi+1}}
=\overline{\pgoth^{-pi}}
\subset_f\cdots
\subset\overline{K}
\leqno{(6)}
$$
(prop.~{9} and~{11}).  Here, $i$ is any integer $>0$, and, as before, an
inclusion of $\Fp$-spaces $E\subset_rE'$ means that $E$ is a codimension-$r$
subspace of $E'$.

There is a canonical isomorphism $\overline\ogoth\to\Fp$ sending $a$ to
$S_{k|\Fp}(\hat b)$, where $\hat b$ is the image in $k/\wp(k)$ of a
representative $b\in\ogoth$ of $a$~; the isomorphism $k/\wp(k)\to\F_p$ is
induced, as before, by the trace map $S_{k|\Fp}$.  

The unramified degree-$p$ extension of $K$ is $K(\wp^{-1}(\ogoth))$
(prop.~{12}).  For an $\F_p$-line\/ $D\neq\overline\ogoth$ in $\overline{K}$
such that $D\subset\overline{\pgoth^{-m}}$ but
$D\not\subset\overline{\pgoth^{-m+1}}$ for some $m>0$, the unique ramification
break of the (cyclic, degree-$p$) extension $K(\wp^{-1}(D))$ occurs at $m$
(prop.~{14}), which is an integer prime to~$p$.

The extension $M=K(\wp^{-1}(K))$ is the maximal elementary abelian
$p$-extension of $K$.  Denote by $G=\Gal(M|K)$ the profinite group of
$K$-automorphism of $M$~; it comes equipped with the ramification filtration
$(G^u)_{u\in[-1,+\infty[}$ in the upper numbering.

We have $G^u=G^1$ for $u\in\;]\!-1,1]$, and, for $u>0$, we have the
``\thinspace orthogonality relation\thinspace''
$(G^u)^\perp=\overline{\pgoth^{-\lceil u\rceil+1}}$ (prop.~{17}), where the
orthogonal is with respect to the Artin-Schreier pairing
$G\times\overline{K}\to\F_p$.  This is the function-field analogue of the
relation $(G^u)^\perp=\bar U_{pe_1-\lceil u\rceil+1}$ (\S1).

The orthogonality relation implies that the upper ramification breaks of $G$
occur precisely at $-1$ and at the integers $>0$ which are prime to~$p$. Given
this, it is tantamount to $K(\wp^{-1}(\pgoth^{-m}))=M^{G^{m+1}}$ for every\/
$m\in\N$ (cor.~{18)}.  This last relation allows us to compute the
discriminant (over~$K$) of these intermediate finite extensions (prop.~{19})~;
the result should be compared with prop.~{3}.

\bigbreak

{\bf 6.  Justifications.}  We now prove the statements of \S5.  As there, $k$
is a finite extension of $\Fp$ of degree~$f$, $K$ is a local field with field
of constants $k$, $\ogoth$ is the ring of integers of $K$ and $\pgoth$ is the
unique maximal ideal of $\ogoth$~; we have $k=\ogoth/\pgoth$.  

We denote by $\wp$ the endomorphism $x\mapsto x^p-x$ of the additive group of
any $\Fp$-algebra, such as $\ogoth$, $K$, $k$.  For any subset $E\subset K$,
denote by $K(\wp^{-1}(E))$ the extension of $K$ obtained by adjoining all
$\alpha$ (in an unspecified algebraic closure of $K$) such that
$\wp(\alpha)\in E$.

Denote by $(\overline{\pgoth^n})_{n\in\Z}$ the filtration on $\overline
K=K/\wp(K)$ induced by the filtration $(\pgoth^n)_{n\in\Z}$ on (the
$\Fp$-space) $K$ (where $\pgoth^0=\ogoth$).

\th LEMMA {8} 
\enonce
We have\/ $\bar\pgoth=\{0\}$.  In fact, $\wp:\pgoth\to\pgoth$ is an
isomorphism. 
\endth
For every $a\in\pgoth$, the reduction $T^p-T$ of the polynomial $T^p-T-a$ has
the $p$ roots $0,1,\ldots,p-1$ making up the subfield $\Fp\subset k$
(``\thinspace Fermat's little theorem\thinspace'').  Hensel's lemma then
implies that there is a unique root $x\in\ogoth$ of $T^p-T-a$ whose reduction
is $0\in k$, so $x\in\pgoth$ is the unique element such that $\wp(x)=a$.  (The
$p$ roots of $T^p-T-a$ are  $x+b$, for $b\in\F_p$)

\medskip

For $n\in\Z$, let $\lambda(n)=\inf(n,pn)$, so that $\lambda(n)=n$ (resp.~$pn$)
if $n\in\N$ (resp.~if $-n\in\N$).  It is clear that
$\wp(\pgoth^n)\subset\pgoth^{\lambda(n)}$ and
$\wp(\pgoth^{n+1})\subset\pgoth^{\lambda(n)+1}$.  There is therefore an
induced map
$\wp_n:\pgoth^n\!/\pgoth^{n+1}\to\pgoth^{\lambda(n)}\!/\pgoth^{\lambda(n)+1}$~;
it is the function-field analogue of the map $(\ )^p:U_n/U_{n+1}\to
U_{\lambda(n)}/U_{\lambda(n)+1}$ in the case of local number fields, where
$\lambda$ was defined only for $n>0$ (as $\inf(pn,n+e)$, $e$ being the
absolute ramification index~; cf.~\citer\locdisc(), Part~III.3).

\th PROPOSITION {9} 
\enonce
The image\/ $\bar\ogoth$ of\/ $\ogoth$ in\/ $\overline K$ is the same as the
cokernel of\/ $\wp_0:\ogoth/\pgoth\to\ogoth/\pgoth$.
\endth
Notice first that $\wp(\ogoth)=\ogoth\cap\wp(K)$.  We have
$\bar\ogoth=\ogoth/\wp(\ogoth)$.  Now, $\pgoth\subset\wp(\ogoth)$ (lemma~{8}),
so $\bar\ogoth$ is also the quotient of $\ogoth/\pgoth$ by
$\wp(\ogoth)/\pgoth=\Im\wp_0$.

To see that $\bar\ogoth$ is canonically isomorphic to $\Fp$, consider the
following commutative diagram
$$
\def\\{\mskip-2\thickmuskip}
\def\droite#1{\\\hfl{#1}{}{8mm}\\}
\diagram{
0&\rightarrow&\F_p&\droite{}&
\ogoth/\pgoth&\droite{\wp_0}
&\ogoth/\pgoth&\droite{}&\bar\ogoth&\rightarrow&0\cr
&&\vfl{=}{}{5mm}&&\vfl{=}{}{5mm}&&\vfl{=}{}{5mm}&&\vfl{}{}{5mm}\cr
0&\rightarrow&\F_p&\droite{}&k&\droite{\wp}
&k&\droite{S_{k|\F_p}}&\F_p&\rightarrow&0\cr
}
$$ 
which is the analogue of the diagram in \citer\locdisc(), Part~III.3 (where
the choice of a primitive $p$-th root of~$1$, or a $(p-1)$-th root of $-p$,
was necessary).   

To bring out the analogy further, consider, for $n\neq0$ in $\Z$, and for
every choice of a uniformiser $\pi$ for $K$, the commutative diagram
$$
\diagram{
\pgoth^n/\pgoth^{n+1}&\droite{\wp_n}&
\pgoth^{\lambda(n)}/\pgoth^{\lambda(n)+1}\cr
\vfl{}{}{5mm}&&\vfl{}{}{5mm}\cr
k&\droite{h}&k,\cr
} 
$$
where $h(x)=-x$ (resp. $x^p$) if $n>0$ (resp.~$n<0$).  The vertical maps are
the isomorphisms induced by the $\ogoth$-bases $\pi^n$, $\pi^{\lambda(n)}$ of
$\pgoth^n$, $\pgoth^{\lambda(n)}$.  In particular, $\wp_n$ is an isomorphism
for $n\neq0$.

\th LEMMA {10}
\enonce
For every integer\/ $n>0$, we have\/ 
$\pgoth^{-n}\cap \wp(K)=\wp(\pgoth^{-\lfloor n/p\rfloor})$.
\endth
Consider $n=pi$ ($i>0$)~; we have to show that
$\pgoth^{-pi}\cap\wp(K)=\wp(\pgoth^{-i})$.  One inclusion follows from
$\wp(\pgoth^{-i})\subset\pgoth^{-pi}$.  For the converse, let $x\in K$ be such
that $\wp(x)\in\pgoth^{-pi}$~; we have to show that $x\in\pgoth^{-i}$.  If
$v(x)\ge0$, then clearly $x\in\pgoth^{-i}$.  Suppose that $v(x)<0$.  Then
$v(\wp(x))=pv(x)$, but by assumption $v(\wp(x))\ge-pi$.  It follows that
$v(x)\ge-i$ and $x\in\pgoth^{-i}$.

Consider next $n=pi+j$ ($i\in\N$, $0<j<p$)~; we have to show that
$\pgoth^{-n}\cap \wp(K)=\wp(\pgoth^{-i})$.  As before, we have
$\wp(\pgoth^{-i})\subset\pgoth^{-pi}\subset\pgoth^{-n}$. If $x\in K$ is such
that $v(x)<0$ and $\wp(x)\in\pgoth^{-n}$, then $v(\wp(x))=pv(x)$, and, as
before, $v(x)\ge-i-(j/p)$.  But $v(x)$ is in $\Z$, so $v(x)\ge-i$ and
$x\in\pgoth^{-i}$.

\th PROPOSITION {11} 
\enonce
Let\/ $m>0$ be an integer.  If\/ $m=pi$ is a multiple of\/ $p$, then\/
$\overline{\pgoth^{-pi+1}}=\overline{\pgoth^{-pi}}$, whereas if\/ $m=pi+j$
($0<j<p$) is prime to\/~$p$, then
$\overline{\pgoth^{-m+1}}\subset\overline{\pgoth^{-m}}$ is a subspace of
codimension\/ $f$ (over\/ $\Fp$).  In particular, the dimension of the
$\Fp$-space\/ $\overline K$ is countably infinite.
\endth
This is the analogue of \citer\locdisc(prop.~42), the major difference being
that for a local number field $F$, the group $F^\times\!/F^{\times p}$ is
finite, and that $\bar U_{pe_1+1}\subset\bar U_{pe_1}$ has codimension~$1$
when $F$ contains a primitive $p$-th root of~$1$.

The proof runs along the same lines.  As there, the source of the dichotomy
between multiples $pi$ of~$p$ and integers $m=pi+j$ ($0<j<p$) prime to~$p$
lies in lemma~{10}, which implies that
$\pgoth^{-pi}\cap\wp(K)=\wp(\pgoth^{-i})$ but
$\pgoth^{-pi+1}\cap\wp(K)=\wp(\pgoth^{-i+1})$, whereas $\pgoth^{-m}\cap
\wp(K)=\wp(\pgoth^{-i})$ and $\pgoth^{-m+1}\cap \wp(K)=\wp(\pgoth^{-i})$.

Consider first multiples of~$p$. In the commutative diagram
$$
\def\\{\mskip-2\thickmuskip}
\def\droite#1{\\\hfl{#1}{}{8mm}\\}
\diagram{
0&\rightarrow&\pgoth^{-i+1}&\droite{}&
\pgoth^{-i}&\droite{}&\pgoth^{-i}/\pgoth^{-i+1}&\rightarrow&0\cr
&&\vfl{\wp}{}{5mm}&&\vfl{\wp}{}{5mm}&&\vfl{\wp_{-i}}{}{5mm}\cr
0&\rightarrow&\pgoth^{-pi+1}&\droite{}&
\pgoth^{-pi}&\droite{}&\pgoth^{-pi}/\pgoth^{-pi+1}&\rightarrow&0\cr
&&\vfl{}{}{5mm}&&\vfl{}{}{5mm}&&\cr
&&\overline{\pgoth^{-pi+1}}&\droite{?}&
\overline{\pgoth^{-pi}}&&\cr
}
$$
the rows and columns are exact and $\wp_{-i}$ is bijective (see above),
because $i\neq0$.  Thus, the arrow marked ``?'' is an isomorphism, by the
snake lemma.

By contrast, for the integer $m$ (prime to $p$), the rows and columns in the
commutative diagram
$$
\def\\{\mskip-2\thickmuskip}
\def\droite#1{\\\hfl{#1}{}{8mm}\\}
\diagram{
0&\rightarrow&\pgoth^{-i}&\droite{=}&\pgoth^{-i}&\droite{}&
 \pgoth^{-i}/\pgoth^{-i}&\droite{=}&0\cr
&&\vfl{\wp}{}{5mm}&&\vfl{\wp}{}{5mm}&&\vfl{0}{}{5mm}\cr
0&\rightarrow&\pgoth^{-m+1}&\droite{}&
    \pgoth^{-m}&\droite{}&\pgoth^{-m}/\pgoth^{-m+1}&\rightarrow&0\cr
&&\vfl{}{}{5mm}&&\vfl{}{}{5mm}&&\cr
&&\overline{\pgoth^{-m+1}}&\droite{}&
\overline{\pgoth^{-m}}&&\cr
}
$$
are as exact as before, but one of the arrows is now~$0$, instead of being an
isomorphism.  Therefore the induced map
$\overline{\pgoth^{-m}}/\overline{\pgoth^{-m+1}}\to\pgoth^{-m}/\pgoth^{-m+1}$
is now an isomorphism, instead of being~$0$.  As the space
$\pgoth^{-m}\!/\pgoth^{-m+1}$ is of dimension $f=[k:\Fp]$, the proof of
prop.~{11} is complete.

{\it Remark}\pointir The same method can be used to determine the filtration on
$\overline{K^\times}=K^\times\!/K^{\times p}$ (which is easily seen to be
infinite, cf.~\citer\locdisc(cor.~21)).  Indeed, we have $U_n^p\subset
U_{pn}$, just as $\wp(\pgoth^{-n})\subset\pgoth^{-pn}$ here.  The result can
be expressed succintly as 
$$
\cdots
\subset_f\bar U_{pi+1}
=\bar U_{pi}\subset_f\cdots
\subset_f\bar U_1
\subset_1\overline{K^\times},
$$
with $\overline{K^\times}/\bar U_1=\Z/p\Z$.  This has the appearance of being
the mirror image of~$(6)$~; the phenomenon will be explained further on.

\th PROPOSITION {12}
\enonce
The unramified degree-$p$ extension of\/ $K$ is\/ $K(\wp^{-1}(\ogoth))$.
\endth

Let $a\in\ogoth$ be such that its image $\bar a$ in $\bar\ogoth$ generates
$\bar\ogoth$.  We have to show that, $\alpha$ being a root of $T^p-T-a$ (in
an algebraic closure of $K$), the extension $K(\alpha)$ is unramified.

But this follows from the fact that the reduction $T^p-T-\hat a\in k[T]$ is an
irreducible polynomial.  Indeed, $S_{k|\Fp}(\hat a)$, being the image of $\bar
a$ under the isomorphism $\bar\ogoth\to\Fp$, generates the latter group, and
hence $\hat a\notin\wp(k)$.

{\it Remark}\pointir Writing $Q=\Gal(K(\alpha)|K)$, we also have the
identification $Q\to\Z/p\Z$ sending $\varphi$ to $1$, where $\varphi$
(``\thinspace Frobenius\thinspace'') is the unique element of $Q$ whose
restriction to $k(\alpha)$ is the $k$-automorphism $\varphi(x)=x^{p^f}$.  In
terms of these identifications, the Artin-Schreier pairing
$Q\times\bar\ogoth\to\Fp$ gets identified with the standard bilinear form 
$(\sigma,c)\mapsto\sigma.c$ from $\Fp\times\Fp$ to $\Fp$.  

This amounts to showing that $\varphi(\beta)-\beta=S_{k|\Fp}(b)$ for every
$b\in k$, where $\beta\in k(\alpha)$ is a root of $T^p-T-b$.  Successively
raising the relation $\beta^p-\beta=b$ to the exponents
$1,p,p^2,\ldots,p^{f-1}$, we get
$$
\beta^p-\beta=b,\ \ 
\beta^{p^2}-\beta^p=b^p,\ \
\ldots,\ \
\beta^{p^f}-\beta^{p^{f-1}}=b^{p^{f-1}}, 
$$
and adding these $f$ equations together gives
$$
\varphi(\beta)-\beta
=\beta^{p^f}-\beta
=b+b^p+\cdots+b^{p^{f-1}}
=S_{k|\Fp}(b),
$$
which was to be proved.

\medbreak

Let us fix some notation.  Let $D\neq\bar\ogoth$ be an $\Fp$-line in
$\overline K$, $m$ the integer such that $D\subset\overline{\pgoth^{-m}}$
but\/ $D\not\subset\overline{\pgoth^{-m+1}}$~; we have seen that $m$ is $>0$
and prime to~$p$ (prop.~{11}).  Fix an element $a\in\pgoth^{-m}$ whose image
generates $D$, let $\alpha$ be a root of $T^p-T-a$ (in an algebraic closure of
$K$), and let $L=K(\alpha)=K(\wp^{-1}(D))$.

Our first task is to find a uniformiser for $L$ (in the analogous case of a
degree-$p$ kummerian extension of local number fields, see \citer\locdisc(),
prop.~61).  We denote the normalised valuations of $K,L$ by $v_K,v_L$~; as the
extension $L|K$ is totally ramified (prop.~{12}) of degree~$p$, we have
$v_L(x)=pv_K(x)$ for every $x\in K$.  Let $\pi$ be any uniformiser of $K$. 

\th PROPOSITION {13}
\enonce
Let\/ $x,y\in\Z$ be such that\/ $-mx+py=1$.  Then\/ $\alpha^x\pi^y$ is a
uniformiser of\/ $L$, and the ring of integers of\/ $L$ is\/
$\ogoth_L=\ogoth[\alpha^x\pi^y]$.
\endth
Notice first that $v_L(\alpha)<0$, for otherwise $\alpha^p-\alpha=a$ would be 
in $\ogoth$.  Therefore $v_L(\alpha^p-\alpha)=v_L(\alpha^p)=pv_L(\alpha)$.  But
we also have $v_L(a)=pv_K(a)=-pm$.  Therefore $v_L(\alpha)=-m$.

It follows that $v_L(\alpha^x\pi^y)=-mx+py=1$, and, because $L|K$ is totally
ramified, that $\ogoth_L=\ogoth[\alpha^x\pi^y]$.

\th PROPOSITION {14}
\enonce
The unique ramification break of the degree-$p$ cyclic extension\/
$L|K$  occurs at\/ $m$.
\endth
Let $H=\Gal(L|K)$ and let $\sigma\in H$ be such that
$\sigma(\alpha)-\alpha=1$~; as $\sigma$ generates $H$, we must show that 
$\sigma\in H_m$ but $\sigma\notin H_{m+1}$.

For this, it is enough \citer\serre(p.~61) to show that
$v_L(\sigma(\varpi)-\varpi)=m+1$ for some uniformiser $\varpi$ of $L$.  We
choose $\varpi=\alpha^x\pi^y$ (prop.~{13}) and compute
$$
{\sigma(\varpi)\over\varpi}
={\sigma(\alpha^x\pi^y)\over\alpha^x\pi^y}
=\left({\sigma(\alpha)\over\alpha}\right)^x
=(1+\alpha^{-1})^x
\equiv1+x\alpha^{-1}\pmod{\varpi^{m+1}},
$$
recalling that $v_L(\alpha^{-1})=m$ and that $x$ is prime to~$p$ (as
$-mx+py=1$).  This shows that $v_L(\sigma(\varpi)-\varpi)=m+1$, hence
$\sigma\in H_m$ but $\sigma\notin H_{m+1}$, hence $H_m=H$ but
$H_{m+1}=\{\Id_L\}$, and the lower (as well as the upper) ramification break
of $H$ occurs at $m$. 

\th COROLLARY {15}
\enonce
The valuation\/ $v_L({\goth D}_{L|K})$ of the different\/ ${\goth D}_{L|K}$
of\/ $L|K$, as well as the  valuation\/ $v_K(d_{L|K})$ of the discriminant,
equals\/ $(p-1)(1+m)$.
\endth
Indeed, $v_L({\goth D}_{L|K})=\sum_{n\in\N}(\Card H_n-1)=(p-1)(1+m)$, and 
$v_K(d_{L|K})$ is the same because $L|K$ is totally ramified.  

The determination of the ramification of a degree-$p$ cyclic extension of $K$
goes back to Hasse \citer\hasse().  I haven't checked if he uses the
uniformiser $\alpha^x\pi^y$.  It wouldn't be surprising if he does, because
$\alpha^x\pi^y$ is just the function-field analogue of $(\xi-\root
l\of\mu)^x\lambda^y$, which can be found in his {\it
  Klassenk{\"o}rperbericht\/}, and even in Hilbert's {\it Zahlbericht} (the
second $\Omega$ in the proof of {\it Satz\/}~148).

{\it Remark}\pointir This allows us --- in principle --- to compute the
discriminant of any finite extension of global function fields.  Briefly, one
reduces first to the local case, then to the galoisian case, then to the case
of a $p$-extension, and finally to the case of a degree-$p$ extension, where
cor.~{15} can be applied.

\medbreak

Now let $M_m=K(\wp^{-1}(\pgoth^{-m}))$ for every $m\in\N$, and
$M=K(\wp^{-1}(K))$, which is the maximal elementary abelian $p$-extension of
$K$.  It is the increasing union
$$
K\subset_1 M_0
\subset_f M_1
\subset_f\cdots
\subset_f M_{pi-1}
= M_{pi}
\subset_f\cdots
\subset M
$$
($i>0$ being arbitrary), where an inclusion of fields $E\subset_rE'$ means
that $E'$ is a degree-$p^r$ extension of $E$.  
\th COROLLARY {16}
\enonce
For every $m\in\N$, the degree of the extension\/ $M_m|K$ is\/ $p^{1+c(m)f}$,
where\/ 
$$
c(m)=m-\left\lfloor{m\over p}\right\rfloor
$$
is the number of integers in\/ $[1,m]$ which are prime to\/ $p$.
\endth
Indeed, the $\Fp$-dimension of $\overline{\pgoth^{-m}}$ is $1+c(m)f$, by
prop.~{11}.  

Define $\displaystyle a(i)=\left\lfloor{i-1\over p-1}\right\rfloor$ as before.
Notice that the strictly increasing map defined by $b^{(i)}=i+a(i)$
establishes a bijection between $[1,c(m)]$ and the set of integers in $[1,m]$
which are prime to $p$.  In other words, the set in question is
$$
b^{(1)}<b^{(2)}<\cdots<b^{(c(m))}.
$$

Put $G_m=\Gal(M_m|M)$ and $G=\Gal(M|K)$~; we are going to think of these
groups as $\Fp$-spaces.  Our next task is to determine the ramification
filtrations on $G_m$ (upper and lower) and on $G$ (upper) in terms of the
Artin-Schreier pairings
$$
G_m\times\overline{\pgoth^{-m}}\to\Fp,\qquad
G\times\overline{K}\to\Fp.
$$
The case $m=0$ is easy~: $M_0|K$ is the unramified degree-$p$ extension
(prop.~{12}), so $G_0^{-1}=G_0$, $G_0^u=\{\Id_{M_0}\}$ for $u>-1$.

\th PROPOSITION {17}
\enonce
We have\/ $G^u=G^1$ for\/ $u\in\;]-1,1]$, and, for\/ $u>0$, 
$$
G^{u\perp}=\overline{\pgoth^{-\lceil u\rceil+1}}
$$
under\/ $G\times\overline{K}\to\Fp$.  The positive ramification breaks in the
filtration on\/ $G$ occur precisely at the integers prime to\/ $p$, namely\/
$b^{(i)}$ ($i>0$).
\endth
Let $u\in\;]-1,1]$.  Notice first that $G^u\neq G$, for otherwise the unique
ramification break of $G/\bar\ogoth^\perp$ would be $>u$, which it is not
(prop.~{12}).  Now let $H$ be a hyperplane containing $G^u$, so that $G/H$ is
cyclic of order~$p$.  As the filtration on $G/H$ is the quotient of the
filtration on $G$, the ramification break of $G/H$ occurs somewhere $<u$
(because $G^u\subset H$).  But the only degree-$p$ cyclic extension of $K$
whose ramification break is $<1$ is $K(\wp^{-1}(\ogoth))$ (prop.~84).  So
$H=\bar\ogoth^\perp$ is the only hyperplane containing $G^u$.  This implies
that $G^u=H=G^1=\bar\ogoth^\perp$.

It remains to show the orthogonality relation
$G^{u\perp}=\overline{\pgoth^{-\lceil u\rceil+1}}$ for $u>1$.  The principle
of the proof is simplicity itself~: two subspaces are the same if they contain
the same lines.  We show that, for a line $D\subset\overline K$, we have
$D\subset G^{u\perp}$ if and only if $D\subset\overline{\pgoth^{-\lceil
    u\rceil+1}}$.

Take a line $D$ and denote by $m$ be the unique ramification break of
$G/D^\perp$.  Then
$$
D\subset G^{u\perp}\Leftrightarrow 
(G/D^\perp)^u=0\Leftrightarrow 
m<u \Leftrightarrow 
m<\lceil u\rceil\Leftrightarrow 
D\subset \overline{\pgoth^{-\lceil u\rceil+1}}.
$$
It now follows from prop.~{11} that the positive ramification breaks of $G$
occur precisely at the integers $b^{(i)}$ ($i>0$) prime to~$p$.  In fact, the
ramification filtration on $G$ looks like
$$
\cdots\subset_f G^{pi+1}=G^{pi}\subset_f\cdots\subset_f G^1=G^0\subset_1G,
$$
($i>0$), where an inclusion $H\subset_r H'$ means that $H$ is a codimension-$r$
subspace (an index-$p^r$ subgroup) of the $\Fp$-space $H'$.

\th COROLLARY {18}
\enonce
For every\/ $m\in\N$, we have\/ $K(\wp^{-1}(\pgoth^{-m}))=M^{G^{m+1}}$.
\endth
Indeed, $(G^{m+1})^\perp=\overline{\pgoth^{-m}}$, by prop.~{17}.

Prop.~{17} allows us to determine the ramification filtration on
$\Gal(L|K)$ for any elementary abelian $p$-extension $L|K$, provided we know
the subspace $\Ker(\overline K\to\overline L)$, where $\overline L=L/\wp(L)$.
Let us do this exercise for $L=K(\wp^{-1}(\pgoth^{-m}))$, in which case
the subspace in question is $\overline{\pgoth^{-m}}$.

\th PROPOSITION {19}
\enonce
The upper ramification breaks of\/ $\Gal(L|K)$ occur at\/ $-1$ and at\/
$b^{(1)}<b^{(2)}<\cdots<b^{(c(m))}$, the $c(m)$ integers in\/ $[1,m]$ which
are prime to\/~$p$.  In the lower numbering, they occur at\/ $-1$ and at
$$
b_{(i)}=(1+q+\cdots+q^{i-1})+(q^{p-1}+\cdots+q^{a(i)(p-1)}), \qquad
i\in[1,c(m)],
$$
where\/ $q=p^f$.  We have \/ $v_{L}({\goth
  D}_{L|K})=(1+b^{(c(m))})q^{c(m)}-(1+b_{(c(m))})$, and\/
$v_{K}(d_{L|K})=pv_{L}({\goth D}_{L|K})$.
\endth
This is very similar to what we saw in prop.~{3} in the case of kummerian
extensions of local number fields.  To compute the valuation $v_{L}({\goth
  D}_{L|K})$ of the different of $L|K$, we appeal to lemma~{2}, noting that the
order of the inertia group $\Gal(L|K)^0$ is $q^{c(m)}$, by cor.~{16}.  The
valuation $v_{K}(d_{L|K})$ of the discriminant is $pv_{L}({\goth D}_{L|K})$
because the residual degree of $L|K$ is $p$.

{\it Remark}\pointir Notice that $L=K(\wp^{-1}(\pgoth^{-m}))$ is the maximal
elementary abelian $p$-extension of $K$ with ramification breaks in $[-1,m]$.
Notice also that the kernel of the projection $G\to\Gal(L|K)$ is $G^{m+1}$
(cor.~{18}), which corresponds to $\bar U_{m+1}$ under the reciprocity
isomorphism $K^\times\!/K^{\times p}\to G$.  But the kernel of
$K^\times\to\Gal(L|K)$ is the group of norms $N_{L|K}(L^\times)$.  It follows
that $N_{L|K}(L^\times)=U_{m+1}K^{\times p}$.

The analogue for a local number field $K$ (containing a primitive $p$-th root
of~$1$) would say that $N_{L|K}(L^\times)=U_{pe_1-m+1}K^{\times p}$ for
$L=K(\!\root p\of{U_m})$, where $m\in[0,pe_1+1]$ and $U_0=K^\times$, which
goes back to \citer\nguyen().

\medbreak

The orthogonality relation of prop.~{17} has other applications as in the case
of local number fields (\S4).  Thus we can determine the ramification
filtration of an arbitrary elementary abelian $p$-extension $L|K$, the
valuation of its discriminant if $L|K$ is finite, and so on.  We can specify
the sequences which can occur as the ramification breaks of some $L|K$, and
determine the possible degrees, and the total number, of such $L|K$.  {\it
  N'insistons pas.}

\bigbreak
\unvbox\bibbox 

\bye